\documentclass[10pt, a4paper]{article}
\usepackage{comment}
\usepackage{graphics}
\usepackage{latexsym}
\usepackage{amssymb}
\usepackage{amsmath}
\usepackage{color}
\usepackage{epsfig}
\usepackage[title]{appendix}
\usepackage{graphicx}
\graphicspath{ {./1D_Uniform/linear_functional_of_a/}{./1D_Uniform/u0_and_u1/}{./1D_Uniform/u0_only/}{./2D_Lognormal/80z/}{./2D_Uniform/} }
\usepackage{caption}
\usepackage{subcaption}

\newcommand{\beq}{\begin{equation}}
\newcommand{\eeq}{\end{equation}}
\newcommand{\beqas}{\begin{eqnarray*}}
\newcommand{\eeqas}{\end{eqnarray*}}
\newcommand{\ep}{\varepsilon}

\newcommand{\ue}{u^{\varepsilon}}

\newcommand{\spa}{{\mathbb R}^d}
\newcommand{\wc}{\rightharpoonup}

\def\IN{\mathbb{N}}

\def\IP{\mathbb{P}}
\def\IR{\mathbb{R}}

\def\IE{\mathbb{E}}

\def\cB{{\cal B}}

\def\cD{{\cal D}}

\def\cG{{\cal G}}

\def\cO{{\cal O}}
\def\cN{{\cal N}}

\def\bproof{{\it Proof}\ \ }
\def\eproof{{\hfill$\Box$}}

\newcommand{\be}{\begin{equation}}
\newcommand{\ee}{\end{equation}}
\newcommand{\bea}{\begin{eqnarray}}
\newcommand{\eea}{\end{eqnarray}}
\newcommand{\beas}{\begin{eqnarray*}}
\newcommand{\eeas}{\end{eqnarray*}}

\newcommand{\bV}{{\bf V}}
\newcommand{\bv}{\boldsymbol{v}}
\newcommand{\bu}{\boldsymbol{u}}

\newcommand{\bw}{\boldsymbol{w}}

\newcommand{\bH}{\boldsymbol{{\cal H}}}

\newcommand{\bg}{\boldsymbol{g}}
\newcommand{\bphi}{\boldsymbol{\phi}}

\newtheorem{theorem}{Theorem}[section]
\newtheorem{lemma}[theorem]{Lemma}
\newtheorem{assumption}[theorem]{Assumption}
\newtheorem{proposition}[theorem]{Proposition}
\newtheorem{definition}[theorem]{Definition}
\newtheorem{remark}[theorem]{Remark}

\definecolor{darkgreen}{rgb}{0,0.7,0}

\numberwithin{equation}{section}
\numberwithin{equation}{section}
\title{\bf Bayesian inverse problems for recovering coefficients of two scale elliptic equations}
\author{
Viet Ha Hoang
and 
Jia Hao Quek\\[10pt]
Division of Mathematical Sciences,\\
School of Physical and Mathematical Sciences,\\
 Nanyang Technological University, Singapore 637371 
}

\date{}
\begin{document}
\maketitle
\begin{abstract}
We consider the Bayesian inverse homogenization problem of recovering the locally periodic two scale coefficient of a two scale elliptic equation, given limited noisy information on the solution. We consider both the uniform and the Gaussian prior probability measures. We use the two scale homogenized equation whose solution contains the solution of the homogenized equation which describes the macroscopic behaviour, and the corrector which encodes the microscopic behaviour. We approximate the posterior probability by a probability measure determined by the solution of the two scale homogenized equation. We show that the Hellinger distance of these measures converges to zero when the microscale converges to zero, and establish an explicit convergence rate when the solution of the two scale homogenized equation is sufficiently regular. Sampling the posterior measure by Markov Chain Monte Carlo (MCMC) method, instead of solving the two scale equation using fine mesh for each proposal with extremely high cost, we can solve the macroscopic two scale homogenized equation. Although this equation is posed in a high dimensional tensorized domain, it can be solved with essentially optimal complexity by the sparse tensor product finite element method, which reduces the computational complexity of the MCMC sampling method substantially. We show numerically that observations on the macrosopic behaviour alone are not sufficient to infer the microstructure. We need also observations on the corrector. Solving the two scale homogenized equation, we get both the solution to the homogenized equation and the corrector. Thus our method is particularly suitable for sampling the posterior measure of two scale coefficients.      
\end{abstract}

\section{Introduction}
%We consider the problem of recovering the microsctructure of a locally periodic two scale problem given limited information on the solution of a two scale elliptic problem posed in the medium. We consider both the cases of uniform prior and log-normal prior.
%
Multiscale problems arise from many important practical and engineering situations such as subsurface flow, reservoir engineering, and composite materials. In many cases, the exact microstructures of the media are not known deterministically. Quantifying the uncertainty in the multiscale media such as finding a description of the permeability of a porous medium is essential for accurate prediction of complex physical processes. The problem of finding  microstructures from limited observations is complex. Observations normally contain random noise. The inverse problems to find the physical properties can be ill-posed, making regularization necessary. Further, solving inverse multiscale problems is highly difficult as it involves solving repeatedly many forward multiscale problems, each of them require large computational resource.

We consider in this paper an inverse problem to find the two scale coefficient of an elliptic equation given limited noisy observations on the solution of the two scale problems. We assume that the coefficient is locally periodic. We follow the framework of Bayesian inverse problems (\cite{KaipioSomersalo}, \cite{AndrewActa}) which assumes that the observation errors follow known probability distributions. Imposing a prior probability measure on the space of locally periodic coefficients, we find the posterior measure which is the conditional probability given the noisy observations. We contribute a method that reduces significantly the complexity of Markov Chain Monte Carlo (MCMC) sampling for the posterior measure in the space of two scale coefficients. 

We consider two types of prior probability measures. In the first type, the two scale coefficients are uniformly coercive and bounded for all the realizations. We assume that they can be written as an expansion of random variables which are uniformly distributed in a compact interval. 
%In the forward uncertainty quantification literature, there is a significant amount of work on random partial differential equations where the coefficients are linear expansions of random variables which are uniformly distributed in a bounded interval (see, e.g., the survey paper \cite{SchwabGittelson}).
 The prior probability space is the probability space of countably infinite sequences of these random variables. We note that the uniform prior probability is considered for Bayesian inversion of single scale elliptic equations in Hoang et al.  \cite{HoangScSt12}. In the second type, the logarithm of the coefficient is a linear expansion of standard normal random vaiables. This is the well known log-gaussian coefficients where the expansion arise from the Kah\`unen-Lo\'eve expansion of the logarithm of the random coefficient given its covariace (see, e.g., \cite{Adler}, \cite{AndrewActa}, \cite{DashtiStuart}). 

As the coefficient of the forward two scale equation is locally periodic, we use two scale convergence (\cite{Nguetseng}, \cite{Allaire}) to study homogenization of the problem. The method produces the two scale homogenized equation for both the solution of the homogenized equation which describes the solution to the two scale equation macroscopically, and the corrector which encodes the microscopic behaviour of the solution to the multiscale equation. These are the first two terms in the two scale asymptotic expansion of the solution to the two scale equation (see \cite{BLP}, \cite{BakhvalovPanasenko}, \cite{JKO}). Solving the two scale homogenized equation, we obtain all the necessary information. The observations obtained from the oscillating solution of the two sale equation can be approximated by the solution of this macroscopic two scale homogenized equation, which leads to an approximation of the posterior measure when the microscopic scale approaches zero. Sampling the posterior by the MCMC method, we can solve this two scale homogenized equation instead of the original two scale equation. Although posed in a high dimensional tensorized domain, the two scale homogenized equation can be solved with an essentially optimal complexity for a prescribed accuracy by the sparse tensor product finite element (FE) method developed by Hoang and Schwab \cite{HSelliptic}. The cost of the MCMC method is far lower than solving the original two scale equation for a large number of samples, using a fine mesh to capture the microscopic scale. Further, we will demonstrate that for recovering the microstructure, it is not sufficient to have only information on the solution of the homogenized problem. We need information on the corrector term. Our method of using the two scale homogenized equation finds all the necessary information with much less complexity. Although we do not consider multilevel MCMC in this paper, the problem is well amenable to this method as developed in \cite{HoangScSt12} which reduces the computation to the optimal level. 

We only consider two scale problems in this paper. However, Bayesian inverse problems for finding coefficients that  depend on multiple separable scales can be solved in the same way, using multiscale convergence and  the essentially optimal method for solving multiscale homogenized equation developed in Hoang and Schwab \cite{HSelliptic}.

Inverse problems for locally periodic two scale problems have been considered before. Nolen et al. \cite{NolenPS} study the effective model and thus the multiscale features of the problem are not recovered. Frederick and Engquist \cite{FrederickEngquist} assume that the microstructure is known, and find a macroscopic quantity such as the volume fraction or the orientation of the inclusions. We contribute in this paper a rigorous theory for recovering the microstructure, and an MCMC sampling method for the posterior measure with significantly reduced complexity. 

For inverse problems of general multiscale media, we mention the work by Efendiev et al. \cite{EfendievHouLuo} (see also the references there in) where a procedure for speeding up the MCMC sampling process  is developed. General multiscale inverse problems using generalized multiscale finite elements is proposed in Chung et al. \cite{ChungEfendiev}.  

The paper is organized as follows. In the next section, we set up the Bayesian inverse problems for two scale equations. We define the prior probability spaces and introduce the Bayes formula for the posterior measure and the well-posedness results of Bayesian inverse problems that we  
will prove later for our two scale setting. We recall the definition of two scale convergence and the results of two scale homogenization for elliptic equations in Section \ref{sec:ts}. We prove that the observations on the oscillating solution of the two scale equation can be approximated by the solution of the two scale homogenized equation, leading to an approximation of the forward functional and the miss-match function. We  study the Bayesian inverse problem with the uniform prior probability in Section \ref{sec:unifprior}. We first establish the existence and well-posedness with respect to the data of the posterior measure.  We show that the posterior measure can be approximated by the corresponding measure determined from the solution of the two scale homogenized equation. In particular, we show that the Hellinger distance of these two mesures converges to zero in the zero limit of the microscale. When the solution of the two scale homogenized equation is sufficiently regular, we establish a rate for this convergence of the Hellinger distance in terms of the microscale. This uses the well-known homgenization rate of convergence. The Bayesian inverse problem with the Gaussian prior probability is studied in Section \ref{sec:gaussprior}.  As in Section \ref{sec:unifprior}, we first prove the existence and well-posedness of the posterior measure. We then prove that the posterior measure can be approximated by the measure determined from the solution of the two scale homogenized equation with respect to  the Hellinger distance, with an explicit bound for the Hellinger distance in terms of the microscopic scale when the solution is sufficiently regular. We review the sparse tensor product FE method for solving the two scale homogenized equation in Section \ref{sec:FE}. This method solves the two scale homogenized equation with essentially equal accuracy as the full tensor product FEs but requires a far less number of degrees of freedom, which is essentially optimal. We define the approximating posterior measure in terms of the FE solutions and prove the convergence of the Hellinger distance of the posterior measure and this FE approximating measure when the meshsize and the microscale converge to zero, with an explicit error bound when the solutions are sufficiently smooth. Numerical examples are presented in Section \ref{sec:numerical}. By considering a one dimensional problem, we show that when we only have observations for the macroscopic behaviour, i.e. when the functions $\ell_i$ in \eqref{eq:Oi} do not depend on $y$, it is not sufficient to recover the microstructure. Indeed, the posterior measure obtained from the MCMC process does not give a clear description of the reference coefficient. Likewise, when the observation is the flux of the two scale equation, which is approximately the flux of the homogenized equation, i.e. we only have an observation for the macroscopic behaviour, the recovery of the reference coefficient by MCMC is poor. However, when the functions $\ell_i$ in \eqref{eq:Oi} depend also on $y$, i.e. we have information on the corrector $u_1$, the MCMC method provides reasonably good recovery of the reference coefficient. We demonstrate this for both the cases of uniform and Gaussian prior probability measures. We conclude the paper with some appendices that contain the long proofs of some results in the previous sections. 

Throughout the paper, by $\nabla$ without indicating the variable, we mean the gradient with respect to $x$ of a function of $x$, and by $\nabla_x$ we mean the partial gradient of a function that depends on $x$ and $y$. Repeated indices indicate summation.  By $c$ we denote a generic constant whose value can change from one appearance to the next. When the constant $c$ depends on the parameter $z$, we write it as $c(z)$. The symbol $\#$ denotes spaces of periodic functions. 

\section{Problem setting}\label{sec:setting}
Let $D\subset\IR^d$ be a bounded domain. Let $Y$ be the unit cube $(0,1)^d\subset\IR^d$. Let $\ep>0$ be a small quantity that represents the microscopic scale of the problem. We consider the prior probability space $(U,\Theta,\rho)$. Let $A:U\to C(\bar D,C_\#(\bar Y))$. We assume that for each $z\in U$, there are constants $c_*(z)>0$ and $c^*(z)>0$ such that 
\be
c_*(z)\le A(z;x,y)\le c^*(z)
\label{eq:coercivebounded}
\ee
for all $x\in D$ and $y\in Y$ \footnote{Indeed we can consider the case where the coefficient $A(z;x,y)$ is a uniformly positive definite symmetric matrix but for simplicity we restrict our consideration to the case where it is a scalar function.}. We define the two scale coefficient as
\[
A^\ep(z;x)=A(z;x,{x\over\ep}).
\]
We denote by $V$ the space $H^1_0(D)$. Let $f\in V'$. 
We consider the forward two scale elliptic problem
\be
-\nabla\cdot(A^\ep(z;x)\nabla u^\ep(z;x))=f(x),\ \ u^\ep=0\mbox{ for }x\in\partial D.
\label{eq:forwardprob}
\ee
Due to condition \eqref{eq:coercivebounded}, problem \eqref{eq:forwardprob} has a unique solution which satisfies
\be
\|u^\ep(z)\|_V\le {c^*(z)\over c_*(z)}\|f\|_{V'}.
\label{eq:bounduep}
\ee
For $N\in \IN$, we consider $N$ functions $\ell_i(x,y)\in C(\bar D,C_\#(\bar Y))^d$ for $i=1,\ldots,N$. We define by 
\be
\cO_i^\ep(z)=\int_D\ell_i(x,{x\over\ep})\cdot\nabla\ue(z;x)dx
\label{eq:Oi}
\ee
which represents a bounded linear map in $H^{-1}(D)$. Let the forward function be
\be
\cG^\ep(z)=(\cO_1^\ep(z),\ldots,\cO^\ep_N(z))\in\IR^N.
\label{eq:Gep}
\ee
We consider a noisy observation $\delta$ of $\cG^\ep(z)$ which is
\[
\delta=\cG^\ep(z)+\nu
\]
where $\nu$ follows a Gaussian distribution $\cN(0,\Sigma)$ in $\IR^N$ where the $N\times N$ covariance matrix $\Sigma$ is positive definite.  We consider the Bayesian inverse problem of determining the posterior probability $\rho^{\delta,\ep}=\IP(z|\delta)$. 

\subsection{Prior probability space}
We consider two types of prior probability measures in this paper. 

\subsubsection{Uniform prior}\label{sec:unif}

We assume that the coefficient $A(z;x,y)$ is represented in the form
\be
A(z;x,y)=\bar A(x,y)+\sum_{i=j}^\infty z_j\psi_j(x,y),
\label{eq:uniform}
\ee
where $\bar A$ and $\psi_j$ belong to $C(\bar D,C_\#(\bar Y))$ with $\inf_{x\in\bar D, y\in\bar Y}\bar A(x,y)>0$. The random variables $z_j$ are uniformly distributed in $[-1,1]$ and are pairwise independent. The prior probability space $U=[-1,1]^\IN$ is equipped with  the $\sigma$ algebra
\[
\Theta=\bigotimes_{j=1}^\infty\cB([-1,1])
\]
where $\cB([-1,1])$ is the Borel $\sigma$ algebra in $[-1,1]$. 
Here $z=(z_1,z_2,\ldots)\in U$. 
The prior probability is
\[
\rho=\bigotimes_{j=1}^\infty{dz_j\over 2}.
\] 
 We assume further that the constants $c_*(z)$ and $c^*(z)$ are uniform with respect to $z$. For this to hold, we make the following assumption.
\begin{assumption}\label{assum:psibound}
There is a constant $\kappa>0$ such that 
\[
\sum_{j=1}^\infty\sup_{x\in D,y\in Y}|\psi_j(x,y)|\le{\kappa\over 1+\kappa}\inf_{x\in D,y\in Y}\bar A(x,y).
\]
\end{assumption}
With this assumption, let
\be
c_*(z)={1\over 1+\kappa}\inf_{x\in D,y\in Y}\bar A(x,y)
\label{eq:clunif}
\ee
and 
\be
c^*(z)=\sup_{x\in D,y\in Y}\bar A(x,y)+ {\kappa\over 1+\kappa}\inf_{x\in D,y\in Y}\bar A(x,y);
\label{eq:cuunif}
\ee
here $c^*(z)$ and $c_*(z)$ do not depend on $z$. 
%Representation \eqref{eq:uniform} can be regarded as an approximation of a \KL expansion of a random field given the covariance, see, e.g. 
%We note that in the forward uncertainty %quantification literature, equations with %parametric coefficients in the form %%\eqref{eq:uniform} have been studied thoroughly %%(see \cite{SchwabGittelson} and the references therein).

Approximation of the solutions of the forward two scale equations with random/parametric coefficients of the form \eqref{eq:uniform} by the generalized polynomial chaos method is studied thoroughly in \cite{HSmultirandom}.
\subsubsection{Gaussian prior}\label{sec:gaussian}
We consider the case where the coefficient $A$ is of the form
\be
A(z;x,y)=A^*(x,y)+\exp(\bar A(x,y)+\sum_{j=1}^\infty z_j\psi_j(x,y)).
\label{eq:gaussian}
\ee
We assume that $z_j$ follow the standard normal distribution ${\cal N}(0,1)$ in $\IR$ and are pairwise independent; here $z=(z_1,z_2,\ldots)\in\IR^N$. 
The function $A^*(x,y)\in C(\bar D,C_\#(\bar Y))$ is non-negative for all $x\in D$ and $y\in Y$; and $\bar A$ and $\psi_j$ belong to $C(\bar D,C_\#(\bar Y))$ for all $j\in\IN$. The functions $\psi_j$ satisfy 
\be
\sum_{j=1}^\infty\sup_{x\in D, y\in Y}|\psi_j(x,y)|<\infty.
\label{eq:finitepsi}
\ee
We equip $\IR^\IN$ with the $\sigma$ algebra
\[
\Theta=\bigotimes_{j=1}^\infty\cB(\IR)
\]
where $\cB(\IR)$ is the Borel $\sigma$ algebra in $\IR$;
and the probability measure
\[
\rho=\bigotimes_{j=1}^\infty \cN(0,1).
\]
We define the space 
\[
U=\{z=(z_1,z_2,\ldots)\in\IR^\IN:\ \sum_{j=1}^\infty |z_j|\sup_{x\in D,y\in Y}|\psi_j(x,y)|<\infty\}.
\]
For conciseness of notation, we denote by $b_j=\sup_{x\in D,y\in Y}|\psi_j(x,y)|$.
With condition \eqref{eq:finitepsi}, we deduce that the set $U$ has  $\rho$ measure 1 (\cite{Yamasaki}, \cite{SchwabGittelson}). We thus define the prior probability space as $U$ with the $\sigma$ algebra $\Theta$ induced in $U$, and the probability $\rho$ restricted to $U$, still denoted as $\rho$. For each $z\in U$, coefficient $A$ is well defined. We also have
\be
A(z;x,y)\ge \inf_{x\in D, y\in Y}A^*(x,y)+\exp(\inf_{x\in D,y\in Y}A(x,y)-\sum_{j=1}^\infty |z_j|b_j):=c_*(z)>0,
\label{eq:cuslogn}
\ee
and 
\be
A(z;x,y)\le \sup_{x\in D, y\in Y}A^*(x,y)+\exp(\sup_{x\in D,y\in Y}A(x,y)+\sum_{j=1}^\infty |z_j|b_j):=c^*(z).
\label{eq:clslogn}
\ee
Problem \eqref{eq:forwardprob} is thus well posed for $z\in U$. We note that if $\inf_{x\in D,y\in Y} A^*(x,y)=0$, $c_*(z)$ is positive for all $z\in U$ but can get arbitrarily close to 0. 
\subsection{Posterior probability measure}
We define the function $\Phi^\ep(z,\delta)$ as 
\[
\Phi^\ep(z,\delta)=\frac12|\delta-\cG^\ep(z)|_\Sigma^2.
\]
In Sections \ref{sec:unifprior} and \ref{sec:gaussprior}, we will show that the posterior probability $\rho^{\delta,\ep}=\IP(z|\delta)$ is determined by the Bayes formula
\be
{d\rho^{\delta,\ep}\over d\rho}\propto\exp(-\Phi^\ep(z,\delta)).
\label{eq:posterior}
\ee
We will show further that $\rho^{\delta,\ep}$ is well posed with respect to $\delta$. In particular, we will show that for $\delta, \delta'\in \IR$ such that $|\delta|\le r$ and $|\delta'|\le r$, 
\be
d_{Hell}(\rho^{\delta,\ep},\rho^{\delta',\ep})\le c(r)|\delta-\delta'|.
\label{eq:localLip}
\ee
Here the Hellinger distance $d_{Hell}$ of two measures $\mu$ and $\mu'$ is defined as (see \cite{AndrewActa})
\[
d_{Hell}(\mu,\mu')=\left(\frac12\int_U\left(\sqrt{d\mu\over d\rho}-\sqrt{d\mu'\over d\rho}\right)^2d\rho\right)^{1/2}.
\]
\section{Two scale convergence}\label{sec:ts}
We review in this section the two scale convergence theory to study homogenization of problem \eqref{eq:forwardprob}. We first recall the definition of two sale convergence which is initiated by Nguetseng \cite{Nguetseng} and developed further by Allaire \cite{Allaire} and Allaire and Briane \cite{AllaireBriane}. 
\begin{definition}
A sequence $\{w^\ep\}_\ep$ in $L^2(D)$ two scale converges to a function $w_0\in L^2(D\times Y)$ if for any functions $\phi\in C(D,C_\#(Y))$
\[
\lim_{\ep\to 0}\int_Dw^\ep(x) \phi(x,{x\over\ep})dx=\int_D\int_Y w_0(x,y)\phi(x,y)dy dx.
\]
\end{definition}
This definition makes sense because of the following result.
\begin{proposition}
From a bounded sequence in $L^2(D)$, we can extract a subsequence which two scale converges.
\end{proposition}
Let $V_1=L^2(D, H^1_\#(Y)/\IR)$. We define by $\bV=V\times V_1$ with the natural norm
\[
\|\bv\|_{\bV}=\|v_0\|_V+\|v_1\|_{V_1},
\]
for $\bv=(v_0,v_1)\in\bV$. 
Using two scale convergence for the two scale elliptic problem \eqref{eq:forwardprob}, we have
\begin{proposition}
The solution $\ue(z;\cdot)$ of \eqref{eq:forwardprob} converges weakly to a function $u_0(z;\cdot)$ in $V$. Further there is a function $u_1(z;\cdot,\cdot)\in V_1$ such that $\nabla\ue$ two scale converges to $\nabla u_0+\nabla_y u_1$.  The function $\bu=(u_0,u_1)\in \bV$ is the unique solution of the problem
\be
B(z;\bu,\bphi)=\int_D\int_Y A(z;x,y)(\nabla u_0(z;x)+\nabla_y u_1(z;x,y))\cdot(\nabla\phi_0(x)+\nabla_y\phi_1(x,y))dy dx=\int_D f(x)\phi_0(x)dx
\label{eq:tsforwardprob}
\ee
$\forall\,\bphi=(\phi_0,\phi_1)\in V\times V_1$. 
\end{proposition}
These results are standard. We refer to \cite{Allaire} for the proof. We then have
\be
\lim_{\ep\to 0}\cO_i^\ep(z)=\int_D\int_Y\ell_i(x,y)\cdot(\nabla u_0(z;x)+\nabla_y u_1(z;x,y))dy dx:=\cO_i^0(z).
\label{eq:Oi0}
\ee
Denoting by
\[
\cG^0(z)=(\cO_1^0(z),\ldots,\cO_N^0(z)),
\]
we have
\be
\lim_{\ep\to 0}\cG^\ep(z)=\cG^0(z).
\label{eq:limcGep}
\ee
Letting
\[
\Phi^0(z,\delta)=\frac12|\delta-\cG^0(z)|^2_\Sigma,
\]
we define the probability measure $\rho^{\delta,0}$ on $U$ as
\be
{d\rho^{\delta,0}\over d\rho}\propto\exp(-\Phi^0(z,\delta)).
\label{eq:rho0}
\ee
We will show that the posterior measure $\rho^{\delta,\ep}$ can be approximated by $\rho^{\delta,0}$ in Sections \ref{sec:unifprior} and \ref{sec:gaussprior} for the uniform prior and the Gaussian prior respectively. We also find a rate of convergence for $d_{Hell}(\rho^{\delta,\ep},\rho^{\delta,0})$ when the solution of \eqref{eq:tsforwardprob} is sufficiently regular.
 
It is well known that the homogenized equation for $u_0$ can be derived from \eqref{eq:tsforwardprob}.
From \eqref{eq:tsforwardprob}, we can write $u_1$ in terms of $u_0$. For $l=1,\ldots,d$, we denote by $w^l(z;x,y)$, as a function of $y\in Y$, the solution of the cell problem
\be
\nabla_y\cdot(A(z;x,y)\nabla_y w^l(z;x,y))=-\nabla_y\cdot (A(z;x,y)e^l),
\label{eq:cell}
\ee
where $w^l(z;x,\cdot)\in H^1_\#(Y)$; $e^l$ is the $l$th unit vector in $\IR^d$ with all the components being zero except the $l$th component which is 1. The symmetric homogenized coefficient $A^0(z;x)$ is determined by
\be
A^0_{kl}(z;x)=\int_YA(z;x,y)(\delta_{kp}+{\partial w^k\over\partial y_p})(\delta_{lp}+{\partial w^l\over\partial y_p}) dy.
\label{eq:A0}
\ee
The function $u_0$ is the solution of the homogenized equation 
\be
-\nabla\cdot(A^0(z;x)\nabla u_0(z;x))=f(x).
\label{eq:homeq}
\ee
The function $u_1$ is determined by
\be
u_1(z;x,y)={\partial u_0\over\partial x_l}(z;x)w^l(z;x,y).
\label{eq:u1}
\ee
%\[
%\lim_{\ep\to 0}d_{Hell}(\rho^{\delta,\ep},\rho^{\delta,0})=0.
%\]
\section{Uniform prior}\label{sec:unifprior}
We consider the case of the uniform prior probability in this section. We first show the existence and well posedness of the posterior probability $\rho^{\delta,\ep}$. We then show the approximation of this measure by the measure $\rho^{\delta,0}$ defined in \eqref{eq:rho0}. 
\subsection{Existence and well-posedness}
We have the following result. 
\begin{proposition}
Under Assumption \ref{assum:psibound}, the posterior probability measure $\rho^{\delta,\ep}=\IP(z|\delta)$ is determined by \eqref{eq:posterior}. Further the posterior measure is well posed, i.e. the local Lipschitz condition \eqref{eq:localLip} holds. 
\end{proposition}
\bproof The proof of this proposition is essentially the same as that of Proposition 3 in \cite{HoangScSt12}. Although Hoang et al. \cite{HoangScSt12} only consider single scale equations but the proof remains valid in our multiscale setting. 
\eproof
\subsection{Approximation by solution of two scale homogenized equation}
We establish the approximation of the posterior measure $\rho^{\delta,\ep}$ by the measure $\rho^{\delta,0}$ defined in \eqref{eq:rho0}.
\begin{theorem}\label{thm:unifconv}
We have
\[
\lim_{\ep\to 0}d_{Hell}(\rho^{\delta,\ep},\rho^{\delta,0})=0.
\]
\end{theorem}
\bproof
Let $Z^\ep(\delta)$ and $Z^0(\delta)$ be the normalizing constants in \eqref{eq:posterior} and \eqref{eq:rho0} respectively, i.e.
\be
Z^\ep(\delta)=\int_U\exp(-\Phi^\ep(z,\delta))d\rho(z),
\label{eq:Zep}
\ee
and 
\be
Z^0(\delta)=\int_U\exp(-\Phi^0(z,\delta))d\rho(z).
\label{eq:Z0}
\ee
We show that $Z^\ep(\delta)$ is uniformly bounded below from 0 for all $\ep$. From \eqref{eq:bounduep}, with $c_*(z)$ and $c^*(z)$ being independent of $z\in U$, we have that $\ue$ is uniformly bounded in $V$. This implies that $\Phi^\ep(z,\delta)\le c(\delta)$ for all $\ep>0$. Thus
$
Z^\ep(\delta)\ge \exp(-c(\delta)).
$
We follow standard procedures of estimating the Hellinger distance of two measures as in \cite{AndrewActa}. We have that
\beqas
2d_{Hell}(\rho^{\ep,\delta},\rho^{0,\delta})^2&=&\int_U\left(Z^\ep(\delta)^{-1/2}\exp(-\frac12\Phi^\ep(z,\delta))-Z^0(\delta)^{-1/2}\exp(-\frac12\Phi^0(z,\delta))\right)^2d\rho(z)\\
&\le&I_1+I_2,
\eeqas
where 
\be
I_1={2\over Z(\delta)}\int_U\left(\exp(-\frac12\Phi^\ep(z,\delta))-\exp(-\frac12\Phi^0(z,\delta))\right)^2d\rho(z),
\label{eq:I1}
\ee
and 
\be
I_2=2|Z^\ep(\delta)^{-1/2}-Z^0(\delta)^{-1/2}|^2\int_U\exp(-\Phi^0(z,\delta))d\rho(z).
\label{eq:I2}
\ee
%Using inequality $|\exp(-x)-\exp(-y)|\le |x-y|$ for $x,y\ge 0$, we have
%\begin{eqnarray}
%|\exp(-\frac12\Phi^\ep(z,\delta))-\exp(-\frac12\Phi^0(z,\delta))|\le c|\Phi^\ep(z,\delta)-\Phi^0(z,\delta)|\nonumber\\
%\le c|\langle\Sigma^{-1/2}(\cG^\ep(z)-\cG^0(z)),\Sigma^{-1/2}(2\delta-\cG^\ep(z)-\cG^0(z))\rangle.
%\label{eq:11}
%\end{eqnarray}
%From \eqref{eq:bounduep}, as $c_*(z)$ and $c^*(z)$ are independent of $z\in U$, we have that $\ue$ is uniformly bounded in $V$. This implies that $|\cG^\ep|$ is uniformly bounded with respect to $\ep$. 
From \eqref{eq:limcGep} we have $\lim_{\ep\to 0}\cG^\ep(z,\delta)=\cG^0(z,\delta)$. Using Lebesgue dominated convergence theorem, we have that
\[
\lim_{\ep\to 0}I_1=0.
\]
We note further that
\beqas
|Z^\ep(\delta)^{-1/2}-Z^0(\delta)^{-1/2}|=\displaystyle{{|Z^0(\delta)-Z^\ep(\delta)|\over Z^\ep(\delta)^{1/2}Z^0(\delta)^{1/2}(Z^0(\delta)^{1/2}+Z^\ep(\delta)^{1/2})}}.
\eeqas
We have 
\[
\lim_{\ep\to 0}|Z^0(\delta)-Z^\ep(\delta)|\le \lim_{\ep\to 0}\int_U\left|\exp(-\Phi^\ep(z,\delta))-\exp(-\Phi^0(z,\delta))\right|d\rho(z)=0.
\]
As $Z^\ep(\delta)>c(\delta)$ uniformly with respect to $\ep$, we deduce that $\lim_{\ep\to 0}I_2=0$. We then get the conclusion. \eproof

With further regularity assumptions on the functions $\psi_j$ in \eqref{eq:uniform}, we get the following convergence rate. We assume:
\begin{assumption}\label{assum:furtherregularitypsiunif}
We assume that the function $\psi_j$ in \eqref{eq:uniform} belong to $C^1(\bar D,C^{1,1}_\#(\bar Y))$ such that
\[
\sum_{j=1}^\infty\|\psi_j\|_{C^1(\bar D,C^{1,1}(\bar Y))}<\infty.
\]
\end{assumption}
\begin{theorem}\label{thm:rateofconvunif} Assume that the domain $D$ is convex, and $f\in L^2(D)$. Under Assumptions \ref{assum:psibound} and \ref{assum:furtherregularitypsiunif}, if $\ell_i\in C^{0,1}(\bar D, C_\#(\bar Y))^d$ then there is a constant $c$ independent of $\ep$ such that
\[
d_{Hell}(\rho^{\delta,\ep},\rho^{\delta,0})\le c\ep^{1/2}.
\]
\end{theorem}
We show this in Appendix \ref{app:a}

\section{Gaussian prior}\label{sec:gaussprior}
We consider the case of the Gaussian prior probability measure in this section. 
\subsection{Existence and well-posedness}
We show the existence and local Lipschitzness with respect to the data of the posterior probability in this section. We have the following result.
\begin{proposition}
With the prior probability space $(U,\Theta,\rho)$ defined in Section \ref{sec:gaussian}, the posterior probability measure $\rho^{\delta,\ep}=\IP(z|\delta)$ is determined by \eqref{eq:posterior}. Further, the local Lipschitz condition \ref{eq:localLip} holds.
\end{proposition}
\bproof
To show that the posterior probability is determined by \eqref{eq:posterior}, we need to show that the forward map $\cG^\ep$ defined in \eqref{eq:Gep} is measurable with respect to the prior $\rho$ (see Cotter et al. \cite{Cotteretal}). The solution $\ue$ of \eqref{eq:forwardprob} is measurable as a map from $(U,\Theta)$ to $(V,\cB(V))$ (see \cite{Gittelson}, \cite{GalvisSarkisSINUM09}, \cite{CharrierSINUM12}). This leads to the measurability of $\cG^\ep$. 

The proof of the local Lipschitz condition \eqref{eq:localLip} follows the same lines of the proof of Proposition 3.3 in \cite{HSmcmclogn}. \eproof

\subsection{Approximation by solution of the two scale homogenized equation}
As for the case of the uniform prior probability measure, we show that the posterior measure $\rho^{\delta,\ep}$ can be approximated by the measure $\rho^{\delta, 0}$ defined in \eqref{eq:rho0}. We first have the following result.
\begin{theorem} \label{thm:gaussconv} For the Gaussian prior measure defined in Section \ref{sec:gaussian}, we have
\[
\lim_{\ep\to 0}d_{Hell}(\rho^{\delta,\ep},\rho^{\delta,0})=0.
\]
\end{theorem}
\bproof
The proof of this theorem is similar to that for Theorem \ref{thm:unifconv}. First we show that the normalizing constant defined in \eqref{eq:Zep} is uniformly bounded from below with respect to all $\ep$, i.e. $Z^\ep(\delta)>c(\delta)>0$ for all $\ep$.. The proof is identical to that for the single macroscopic scale problem in \cite{HSmcmclogn}.

 Indeed, from \eqref{eq:bounduep}, \eqref{eq:cuslogn} and \eqref{eq:clslogn}, we have
\be
\|\ue(z)\|_V\le c\exp(2\sum_{j=1}^\infty|z_j|b_j)).
\label{eq:bounduegauss}
\ee
Thus $|\Phi^\ep(z,\delta)|\le c(1+\exp(4\sum_{j=1}^\infty|z_j|b_j))$. This implies that $\IE^\rho(\Phi^\ep(z,\delta))<\Lambda$ for a positive constant $\Lambda$ (see Lemma \ref{lem:B100} in Appendix \ref{app:b}). Fixing a constant $C$, the measure of the set $\{z\in U: \Phi^\ep(z,\delta)>C\}$ is not more than $\Lambda/C$. Choosing $C$ sufficiently large so that $1-\Lambda/C>0$, we have
\[
Z^\ep(\delta)\ge \exp(-C)(1-\Lambda/C).
\]
We show that $\lim_{\ep\to 0}I_1=0$ and $\lim_{\ep\to 0}I_2=0$ where $I_1$ and $I_2$ are defined in \eqref{eq:I1} and \eqref{eq:I2} respectively. From \eqref{eq:limcGep}, we have that
\[
\lim_{\ep\to 0}\left|\exp(-\frac12\Phi^\ep(z,\delta))-\exp(-\frac12\Phi^0(z,\delta))\right|=0.
\]
From Lebesgue dominated convergence theorem, we deduce that $\lim_{\ep\to 0}I_1=0$. Similarly, $\lim_{\ep\to 0}I_2=0$. \eproof

With regularity conditions, we show a rate for the convergence in the previous theorem. We make the following assumption.
\begin{assumption}\label{assum:furtherregularitypsigauss}
We assume that the functions $A^*$, $\bar A$ and  $\psi_j$ in \eqref{eq:gaussian} belong to  $C^1(\bar D,C^{1,1}_\#(\bar Y))$ such that
\[
\sum_{j=1}^\infty\|\psi_j\|_{C^1(\bar D,C^{1,1}(\bar Y))}<\infty.
\]
\end{assumption}
Let $\bar b_j=\|\psi_j\|_{C^1(\bar D, C^{1,1}(\bar Y))}$. Let $\bar U\subset U$ be the set of $z=(z_1,z_2,\ldots)\in U$ such that $\sum_{j=1}^\infty |z_j|\bar b_j$ is finite. We have that $\rho(\bar U)=1$. For $z\in\bar U$, the coefficient $A$ defined in \eqref{eq:gaussian} belongs to $C^1(\bar D, C^1_\#(\bar Y))$ with
\[
\|A(z)\|_{C^1(\bar C,C^{1,1}(\bar Y))}\le c(1+\exp(c\sum_{j=1}^\infty |z_j|\bar b_j))
\]
for $z\in\bar U$. 
 We have the following result.
\begin{theorem}\label{thm:rateofconvgauss}
Assume that $D$ is convex, and $f\in L^2(D)$. Under Assumption \ref{assum:furtherregularitypsigauss}, if $\ell_i\in C^{0,1}(\bar D, C_\#(\bar Y))^d$ then  
\[
d_{Hell}(\rho^{\delta,\ep},\rho^{\delta,0})\le c\ep^{1/2}.
\]
\end{theorem}
We prove this theorem in Appendix \ref{app:b}.

\section{Finite Element approximation of the posterior measure}%method for the two scale homogenied forward problem \eqref{eq:tsforwardprob}}
\label{sec:FE}
Sampling the probability $\rho^{\delta,0}$, we get an approximation for the posterior probability $\rho^{\delta,\ep}$. The advantage of sampling $\rho^{\delta,0}$ is that we only need to solve equation \eqref{eq:tsforwardprob} which only involves the macroscopic scale. Although this problem is posed in the high dimensional tensor product domain $D\times Y$, the sparse tensor product FE method is capable of solving this problem with an essentially optimal complexity to obtain an essentially equal accuracy as that of the full tensor product FE spaces. We review in this section the sparse tensor product FE method originally developed in \cite{HSelliptic}.
\subsection{Hierarchical finite element spaces} \label{subsec:H}
Let $D$ be a polyhedron. We consider in  $D$ a hierarchy $\{{\cal T}^l\}$ of simplices which are obtained recursively.  The set of simplices ${\cal T}^l$ of mesh size $h_l=O(2^{-l})$ is obtained by dividing each simplex in ${\cal T}^{l-1}$ into 4 congruent triangles for $d=2$ or 8 tedrahedra for $d=3$. The cube $Y$ is divided into sets of simplices ${\cal T}_\#^l$ which are periodically distributed in a similar manner. We define the following finite element spaces  
\beqas
V^l=\{w\in H^1(D): w|_K\in {\cal P}_1(K)\ \forall K\in {\cal T}^l\},\\
V^l_0=\{w\in H^1_0(D): w|_K\in {\cal P}_1(K)\ \forall K\in {\cal T}^l\},\\
V^l_\#=\{w\in H^1_\#(Y): w|_K\in {\cal P}_1(K)\ \forall K\in {\cal T}^l_\#\},
\eeqas
where ${\cal P}_1(K)$ denotes the set of linear polynomials in $K$. The following approximation properties hold (see, e.g. Ciarlet \cite{Ciarlet}).
%
%
%
%We consider  hierarchies of FE subspaces of $H^1_0(D)$, $H^1(D)$ and $H^1_\#(Y)$ 
%\beqas
%V_0^0\subset V_0^1\subset\ldots\subset V_0^l\subset\ldots\subset H^1_0(D),\\
%V^0\subset V^1\subset\ldots\subset V^l\subset\ldots\subset H^1(D),\\
%V_\#^0\subset V_\#^1\subset\ldots\subset V_\#^l\subset\ldots\subset H^1_\#(Y)
%\eeqas
%so that $\bigcup_{l=1}^\infty V_0^l$ is dense in $H^1_0(D)$, $\bigcup_{l=1}^\infty V^l$ is dense in $L^2(D)$ and $\bigcup_{l=1}^\infty V_\#^l$ is dense in $H^1_\#(Y)$. These spaces consist of piecewise linear (or linear with respect to each coordinates in the case of rectangle elements) functions  with the mesh size $h_l=O(2^{-l})$ so that the following approximation properties hold (see Ciarlet \cite{Ciarlet}. For $s\ge 0$
\be
\begin{array}{lr}
\displaystyle
\inf_{v\in V^l}\|w-v\|_{H^1(D)}\le ch_l^{\min\{1,s\}}\|w\|_{H^{1+s}(D)},\\
\displaystyle
\inf_{v\in V^l}\|w-v\|_{L^2(D)}\le ch_l^{1+\min\{1,s\}}\|w\|_{H^{1+s}(D)},
\end{array}
\label{eq:e1}
\ee
for all $w\in H^{1+s}(D)$, and
\be
\inf_{v\in V_0^l}\|w-v\|_{H^1(D)}\le ch_l^{\min\{1,s\}}\|w\|_{H^{1+s}(D)}
\ee
for $w\in H^{1+s}(D)\bigcap V$.  For periodic functions in $Y$ we have
\be
\begin{array}{lr}
\displaystyle
\inf_{v\in V^l_\#}\|w-v\|_{H^1(Y)}\le ch_l^{\min\{1,s\}}\|w\|_{H^{1+s}(Y)},\\
\displaystyle
\inf_{v\in V^l_\#}\|w-v\|_{L^2(Y)}\le ch_l^{1+\min\{1,s\}}\|w\|_{H^{1+s}(Y)},
\end{array}
\ee
for all periodic functions $w\in H^{1+s}_\#(Y)$. 

%{\bf Example}\ We divide $D$ into a hierarchy $\{{\cal T}^l\}$ of simplices recursively.  The set of simplices ${\cal T}^l$ of mesh size $h_l=O(2^{-l})$ is obtained by dividing each simplex in ${\cal T}^{l-1}$ into 4 congruent triangles for $d=2$ or 8 tedrahedra for $d=3$. The cube $Y$ is divided into sets of simplices ${\cal T}_\#^l$ which are periodically distributed in a similar manner. We define the following finite element spaces  
%
%\beqas
%V^l=\{w\in H^1(D): w|_K\in {\cal P}_1(K)\ \forall K\in {\cal T}^l\},\\
%V^l_0=\{w\in H^1_0(D): w|_K\in {\cal P}_1(K)\ \forall K\in {\cal T}^l\},\\
%V^l_\#=\{w\in H^1_\#(Y): w|_K\in {\cal P}_1(K)\ \forall K\in {\cal T}^l_\#\},
%\eeqas
%where ${\cal P}_1(K)$ denotes the set of linear polynomials in $K$. The approximation properties above hold for these spaces (see, e.g. Ciarlet \cite{Ciarlet}).

Similarly, if $D$ is a union of squares ($d=2$) or cubes ($d=3$), we first divide $D$ into subsquares. The simplices in ${\cal T}^l$ are obtained by dividing each simplex in ${\cal T}^{l-1}$ into 4 squares for $d=2$ or 8 cubes for $d=3$. The simplices in $Y$ are defined similarly. The finite element spaces consist of functions that are linear with respect to each component of $x$ or $y$ in each simplex, i.e. they belong to the space ${\cal Q}_1(K)$ for $K\in {\cal T}^l$ and $K\in {\cal T}^l_\#$.    
\subsection{Full tensor product finite elements for problem \eqref{eq:tsforwardprob}}\label{sec:full}
As $u_{1}\in L^2(D,H^1_\#(Y))\cong L^2(D)\otimes H^1_\#(Y)$, we choose the finite element subspace $V_1^L \subset V_1$ for approximating $u_{1}$ as
\be
\label{eq:V1l}
V_1^L=V^L\otimes V^L_\#.
\ee
Let
\be
{\bV}^L = V^L_0 \times V^L_1.
\ee
We consider the full tensor product approximating problem: Find $\bu^L(z)=(u_0^L(z),u_1^L(z))\in \bV^L$ such that
\be
B(z;\bu^L(z),\bv^L)=\int_Dfv_0^Ldx\ \ \ \forall\,\bv^L=(v_0^L,v_1^L)\in \bV^L.
\label{eq:fulltensor}
\ee

To get an explicit error estimate for the finite element approximation, we define the regularity space ${\bH}$ as follows. 

Let ${\cal H}$ be the space of functions $w(x,y)\in L^2(D, H^2_\#(Y))$ such that $w\in L^2(Y,H^1(D))$. The space ${\cal H}$ is equipped with the norm
\beqas
\|w\|_{\cal H}=\|w\|_{L^2(D, H^2_\#(Y))}+\|w\|_{L^2(Y,H^1(D))}.
\eeqas
We then have the following approximating property whose proof can be found in \cite{BungartzGriebel}, \cite{HSelliptic} 
\begin{lemma}
\label{lem:Hiregularity}
 For $w\in{\cal H}$, 
\[
\inf_{v\in V_1^L}\|w-v\|_{L^2(D,H^1(Y))}\le ch_L\|w\|_{\cal H}.
\]
\end{lemma}
We define the space $\bH$ by
\be
\bH=H^2(D)\times {\cal H}, 
\label{eq:bH}
\ee
which is equipped with the norm
\[
\|\bw\|_{\bH}=\|w_{0}\|_{H^2(D)}+\|w_1\|_{\cal H},
\]
where $\bw=(w_0,w_1)\in\bH$.

From this we deduce the following rate of convergence for the full tensor approximating problem \eqref{eq:fulltensor}.
\begin{proposition}\label{prop:fullerror}
Assume that the solution $\bu(z)=(u_0(z),u_1(z))$ of problem \eqref{eq:tsforwardprob} belongs to $\bH$, then for the FE approximating problem \eqref{eq:fulltensor}, we have
\[
\|\bu(z)-\bu^L(z)\|_{\bV}\le c{c^*(z)\over c_*(z)}h_L\|\bu(z)\|_{\bH}.
\]
\end{proposition}
%{\it Proof}\ \ Using Cea's lemma, we get
%\[
%\|\bu-\bu^L\|_{\bV}\le c\inf_{\bv^L\in\bV^L}\|%\bu-\bv^L\|_{\bV}.
%\]
%We then get the conclusion from Lemma \ref{lem:Hiregularity}.  \hfill$\Box$
The proof uses Cea's lemma and Lemma \ref{lem:Hiregularity}. 
 
\subsection{Sparse tensor product finite elements for problem \eqref{eq:tsforwardprob}}\label{sec:sparse}
The dimension of the full tensor product space $\bV^L$ is $O(2^{dL})$ which is prohibitively large when $L$ is  large. %for solving problem \eqref{eq:fulltensor}. 
We develop the sparse tensor product finite element spaces with an essentially optimal  dimension but produce essentially equal accuracy as for the full tensor product FE spaces.  
%%%%%%%%%%%%%%%%%%%%%%%%%%%%%%%%%%%%%%%%%%%%%%%%%%%%%%%%%%%%%%%%%%%%%%%
%We first define the following orthogonal projections in the norm of $L^2(D)$, $L^2(Y)$ and $H^1_{\#}(Y)$
%\beqas
%P^{l0}:L^2(D)\to V^l,\ \ P_\#^{l0}:L^2(Y)\to V_\#^l,\ \ P^{l1}_\#:H^1_\#(Y)\to V_\#^l.
%\eeqas
%With the convention that $P^{-10}=0$, $P^{-10}_\#=0$ and $P^{-11}_\#=0$, the increment spaces are defined as
%\[
%W^l=(P^{l0}-P^{(l-1)0})V^l,\ \ W_\#^{l0}=(P_\#^{l0}-P_\#^{(l-1)0})V_\#^l,\ \ W_\#^{l1}=(P_\#^{l1}-P_\#^{(l-1)1})V_\#^l.
%\]
%We then have
%\[
%V^l=\bigoplus_{0\le i\le l}W^i,\ \ V^l_\#=\bigoplus_{0\le i\le l}W_\#^{i0}=\bigoplus_{0\le i\le l}W_\#^{i1}.
%\]
%%%%%%%%%%%%%%%%%%%%%%%%%%%%%%%%%%%%%%%%%%%%%%%%%%%%%%%%%%%%%%%%%%%%%%%%
We assume that  for each $l\ge 1$ there is a linear space $W^l\subset V^l$ that is linearly independent of $V^{l-1}$ so that $V^l$ is the linear span of $W^l$ and $V^{l-1}$. We denote this as $V^l=V^{l-1}\bigoplus W^l$
% (the notation $\bigoplus$ here does not necessarily indicate an orthogonal summation).   
 with ${\rm dim}V^l={\rm dim}W^l+{\rm dim}V^{l-1}$. Here $W^0=V^0$. Let $\{\psi^l_j\}$ where $j\in I_l\subset\IN$ be a linear basis of $W^l$. Then  $\{\psi^{l'}_j\}$ for $0\le l'\le l$ and $j\in I_{l'}$ form a basis for $V^l$. We assume further that this is a Riesz wavelet basis in $L^2(D)$, i.e. there are constants $c_1>0$ and $c_2>0$ such that for all $w=\sum_{l=0}^\infty\sum_{j\in I_l}w^l_j\psi^l_j$, we have the norm equivalence
\be
c_1\sum_{l=0}^\infty\sum_{j\in I_l}(w^l_j)^2\le \|w\|_{L^2(D)}^2\le c_2\sum_{l=0}^\infty\sum_{j\in I_l}(w^l_j)^2.
\label{eq:normequiv}
\ee
For $H^1_\#(Y)/\IR$, we suppose that there are spaces  $W_\#^{l1}={\rm span}\{\psi_{1j}^l:\ \ j\in I_{1l}\subset\IN\}$ with $W_\#^{01}=V_\#^0$ such that $V^l_\#={\rm span}\{\psi_{1j}^{l'}:\ \ l'=0,\ldots,l\mbox{ and } j\in I_{1{l'}}\}$. We further have the norm equivalence
\be
c_5\sum_{l=0}^\infty\sum_{j\in I_{1l}}(w^l_j)^2\le \|w\|_{H^1(Y)/\IR}^2\le c_6\sum_{l=0}^\infty\sum_{j\in I_{1l}}(w^l_j)^2
\label{eq:normequiv2}
\ee
for $w\in H^1_\#(Y)/\IR$ where $c_5$ and $c_6$ are independent of $w$. %We define
%\be
%Q^l_1w=\sum_{j\in I_{1l}}w^l_j\psi^l_{1j}
%\label{eq:ql1}
%\ee
%with the approximation property 
%\be
%\|Q^l_1w\|_{H^1_\#(Y)/\IR}\le ch_l\|w\|_{H^2(Y)}
%\label{eq:ee3}
%\ee
%when $w\in H^2(Y)$. 
From 
\[
V^l=\bigoplus_{0\le l'\le l}W^{l'},\ \ V^l_\#%=\bigoplus_{0\le l'\le l}W^{l'0}_\#
=\bigoplus_{0\le l'\le l}W^{l'1}_\#,
\]
the full tensor product space $ V_1^L$ in \eqref{eq:V1l}  can be written as
\[
V_1^L=\bigoplus_{{0\le l_0,l_1\le L}}W^{l_0}%\otimes W^{l_10}_\#\otimes\ldots\otimes W^{l_{m-1}0}_\#
\otimes W_\#^{l_11}.
\]
We define the sparse tensor product space $\hat V^L_1\subset V^L_1$ as
\[
\hat V_1^L=\bigoplus_{0\le l_0+l_1\le L}W^{l_0}%\otimes W^{l_10}_\#\otimes\ldots\otimes W^{l_{m-1}0}_\#
\otimes W_\#^{l_11},
\]
and the finite element space
\[
\hat\bV^L=V^L_0%\times (\hat{V}_1^L)^d\times\ldots
\times \hat{V}_1^L.
%\{(\hat u_{0,i}^L,\{\hat u_{p,i}^L\}):\hat u_{0,i}^L,\hat u_{p,i}^L\in \hat V_p^L\}.
\]
%for approximating the solution $\bu$ of \eqref{eq:tsforwardprob}. 
We consider the sparse tensor product FE approximating problem: Find $\hat\bu^L(z) \in \hat\bV^L$ such that
\begin{equation}
B(z;\hat\bu^L(z),\hat\bv^L)= \int_Df(x)\hat v_0^L(x)dx
\label{eq:sparsetensor}
\end{equation}
for all $\hat\bv^L=(\hat v_0^L,\hat v_1^L)\in \hat{\bV}^L$.
%By \eqref{eq:coercivebounded}, we deduce from Lax-Milgram lemma and Korn inequality that this problem has a unique solution. 
%We consider the semi-discretized finite element approximation problem \eqref{eq:discreteproblem} and the fully discretized problem \eqref{eq:discretetimeproblem} with the spaces $V_i^L$ being replaced by $\hat V_i^L$ and $\bV^L$ being replaced by $\hat\bV^L$. 
To get a FE rate of convergence for problem \eqref{eq:sparsetensor}, we define the regularity spaces $\hat{\cal H}$ of functions $w(x, y)$ that are $Y$-periodic with respect to $y$  such that for all $\alpha_0 ,\alpha_1 \in \IN_0^d$ with $|\alpha_0 | \le 1$ and $|\alpha_1 |\le 2$, 
\[
{\partial^{ |\alpha_0 |+|\alpha_1|} w\over\partial^{\alpha_0} x\partial^{\alpha_1} y }\in L^2 (D \times Y),
\]
i.e. $w\in H^1(D,H^2_\#(Y))$. 
We then equip the space $\hat{\cal H}$ with the norm
\[
\|\bw\|_{\hat{\cal H}}=\sum_{|\alpha_1|\le 2, |\alpha_0|\le 1}\left\|{\partial^{ |\alpha_0 |+|\alpha_1 |} w\over\partial^{\alpha_0} x\partial^{\alpha_1} y}\right\|_{ L^2 (D \times Y)}.
\]
We define the regularity space
\[
\hat\bH=H^2(D)\times\hat{\cal H}
\]
with the norm
\[
\|\bw\|_{\hat\bH}=\|w_0\|_{H^2(D)}+\|w_1\|_{\hat{\cal H}}
\]
for $w=(w_0,w_1)\in V\times V_1$. 
For functions in $\hat{\cal H}$,  we have the following estimate
% This type of results for sparse tensor product approximations is well known 
(see  e.g., \cite{BungartzGriebel}, \cite{vonPS}, for a proof). %For completeness, we present the proof here.
\begin{lemma} \label{lem:sparseerror}For $w\in\hat{\cal H}$:
\[
\inf_{v\in \hat V_1^L}\|w-v\|_{V_1}\le cL^{1/2}h_L\|w\|_{\hat{\cal H}}.
\]
\end{lemma}

We now present some examples of the wavelet basis functions that satisfy the norm equivalence above.
%\eqref{eq:normequiv}, \eqref{eq:normequiv1} and \eqref{eq:normequiv2}. 
For constructions of wavelet basis functions we refer to references such as  \cite{Cohenwavelet} and \cite{Dahmenwavelet}. 

{\bf Example}
(i)  A hierarchical basis for $L^2(0,1)$ can be constructed as follows.  We first take three following piecewise linear functions as the basis for level $l=0$: $\psi^{0}_1$ obtains values $(1,0)$ at $(0,1/2)$ and is 0 in $(1/2,1)$, $\psi^{0}_2$ is continuous piecewise linear and obtains values $(0, 1, 0)$ at $(0, 1/2, 1)$, and $\psi^{0}_3$ obtains values $(0,1)$ at $(1/2, 1)$ and is 0 in $(0,1/2)$. The basis functions for other levels are constructed from the wavelet function $\psi$ that takes values $(0,-1,2,-1,0)$ at $(0,1/2,1,3/2,2)$, the left boundary function $\psi^{left}$ taking values $(-2,2,-1,0)$ at $(0,1/2,1,3/2)$, and the right boundary function $\psi^{right}$ taking values $(0, -1,2,-2)$ at $(1/2,1,3/2,2)$. For levels $l\geq 1$, $I_l=\{1,2,\ldots,2^l\}$. The wavelet basis functions are defined as $\psi^{l}_1(x) = 2^{-l/2}\psi^{left}(2^l x)$, $\psi^{l}_k(x)=2^{-l/2}\psi(2^l x - k + 3/2)$ for $k = 2, \cdots, 2^l-1$ and $\psi^{l}_{2^l} = \psi^{right}(2^l x - 2^l+2)$. %This base satisfies Assumption \ref{FEassumption} (i).\\
%(ii) For space $H^1_0(D)$, a simple hierarchical basis function that satisfies Assumption \ref{FEassumption} can be constructed from  the hat function that is piecewise linear and  obtains the values $(0,1,0)$ at $(0,1/2,1)$.
%At level $j$, $I^j=\{1,2,\ldots,2^l\}$ and  $\phi^{jk}(y)=2^{-j/2}\phi(2^jy-i+1)$. 

(ii) For $Y = (0,1)$, a hierarchical basis  for $H^1_{\#}(Y)/\IR$ can be constructed from those in  (i). For level 0, we exclude $\psi^{0}_1$, $\psi^{0}_3$. At other levels, the  functions  $\psi^{left}$ and $\psi^{right}$ are replaced by the piecewise linear functions that take values $(0,2, -1, 0)$ at $(0,1/2,1,3/2)$ and values $(0, -1,2,0)$ at $(1/2,1, 3/2 ,2)$ respectively.

%Another hierarchical basis  can be constructed from the simple hat function that obtains the values $(0,1,0)$ at $(0,1/2,1)$. The wavelet functions $\phi^{jk}(y)=2^{-j/2}\phi(2^jy-i+1)$.

For the $d$ dimensional cube $(0,1)^d$, the basis functions can be constructed by taking the tensor products of the basis functions in $(0,1)$. They satisfy the norm equivalence after appropriate scaling, see \cite{GO95}.

We then have the following error estimate for sparse tensor product FE  approximation.
\begin{proposition}\label{prop:sparseapprox}
Assume that the solution $\bu$ of \eqref{eq:tsforwardprob} belongs to $\hat\bH$. Then the solution $\hat\bu^L$ of the finite element approximating problem \eqref{eq:sparsetensor} satisfies
\be
\|\bu(z)-\hat\bu^L(z)\|_{\bV}\le c{c^*(z)\over c_*(z)}L^{1/2}h_L\|\bu(z)\|_{\hat\bH}.
\label{eq:sparseapprox}
\ee
\end{proposition}
The proof of this proposition uses Cea's lemma and Lemma \ref{lem:sparseerror}.
\subsection{FE approximation of the posterior measure}
We define an approximation to the posterior measure using the FE solution $\hat\bu^L(z)=(\hat u_0^L(z),\hat u_1^L(z))$ in \eqref{eq:sparsetensor}. We define for $i=1,\ldots,N$
\[
\cO_i^{0,L}(z)=\int_D\int_Y\ell_i(x,y)\cdot(\nabla\hat u_0^L+\nabla_y\hat u_1^L)dydx;
\]
and 
\[
\cG^{0,L}(z)=(\cO_1^{0,L}(z),\ldots,\cO_N^{0,L}(z)).
\]
The function $\Phi^{0,L}$ which approximates $\Phi^0$ is defined as
\[
\Phi^{0,L}(z,\delta)=\frac12|\delta-\cG^{0,L}(z)|^2_\Sigma.
\]
We define the measure $\rho^{\delta,0,L}$ as
\[
{d\rho^{\delta,0,L}\over d\rho}\propto\exp(-\Phi^{0,L}(z,\delta)).
\]
We then have the following result.
\begin{theorem}\label{thm:FEconv}
For both the cases of uniform prior measure \eqref{eq:uniform} and Gaussian prior measure \eqref{eq:gaussian}, we have the following approximation
\[
\lim_{\ep\to 0, L\to\infty} d_{Hell}(\rho^{\delta,\ep},\rho^{\delta,0,L})=0.
\]
\end{theorem}
\bproof As $\lim_{L\to\infty}\Phi^{0,L}(\delta,z)=\Phi^0(\delta,z)$, by arguing as in the proofs of Theorems \ref{thm:unifconv} and \ref{thm:gaussconv} we have
\[
\lim_{L\to\infty}d_{Hell}(\rho^{\delta,0},\rho^{\delta,0,L})=0.
\]
We then get the conclusion from Theorems \ref{thm:unifconv} and \ref{thm:gaussconv}. 
\eproof

With regularity assumptions, we get an explicit rate of convergence for the convergence in Theorem \ref{thm:FEconv}. 
\begin{theorem}
Assume that $D$ is a convex polygon and $f\in L^2(D)$. We assume Assumption \ref{assum:furtherregularitypsiunif} for the case of uniform prior \eqref{eq:uniform} and Assumption \ref{assum:furtherregularitypsigauss} for the case of Gaussian prior \eqref{eq:gaussian}. Then
\[
d_{Hell}(\rho^{\delta,\ep},\rho^{\delta,0,L})\le c(\ep^{1/2}+L^{1/2}2^{-L}).
\]
\end{theorem}
\bproof For the case of uniform prior probability \eqref{eq:uniform}, under Assumption \ref{assum:furtherregularitypsiunif} we have that $\hat\bu(z)=(\hat u_0(z),\hat u_1(z))$ is uniformly bounded in $\hat{\bH}$ (see \cite{HSmultirandom} Proposition 4.2). Thus
\[
|\cG^0(z)-\cG^L(z)|\le c\|\bu(z)-\hat\bu^L(z)\|_{\bV}\le cL^{1/2}2^{-L}\ \ \forall\,z\in\bar U,
\]
where $|\cdot|$ denotes the Euclidean norm in $\IR^N$.
% We note that
%\beqas
%|\exp(-\frac12\Phi^0(z,\delta)-\exp(-\frac12\Phi^0(z;%\delta)|\le c|\Phi^0(z,\delta)-\Phi^{0,L}(z,\delta)|\\
%\le c|\langle\Sigma^{-1/2}(\cG^0(z)-\cG^{0,L}(z)),\Sigma^{-1/2}(2\delta-\cG^0(z)-\cG^{0,L}(z))\rangle\le cL^{1/2}2^{-L}.
%\eeqas 
By a similar procedure as in the proof of  Theorem \ref{thm:rateofconvunif}, using $\rho(\bar U)=1$, we show that
\[
d_{Hell}(\rho^{\delta,0},\rho^{\delta,0,L})\le cL^{1/2}2^{-L}.
\]
From this and Theorem \ref{thm:rateofconvunif} we get the conclusion. 

For the Gaussian prior in \eqref{eq:gaussian}, we have that 
\[
|\cG^0(z)-\cG^L(z)|\le c\|\bu(z)-\hat\bu^L(z)\|_{\bV}\le c\|\bu(z)\|_{\hat{\bH}}L^{1/2}2^{-L},
\]
From Lemmas \ref{lem:wlnorm} and \ref{lem:H2normu0}, and from \eqref{eq:u1}, we have that 
\[
\|\bu(z)\|_{\hat{\bH}}\le c(1+\exp(c\sum_{j=1}^\infty|z_j|\bar b_j))
\]
for $z\in\bar U$. By an identical proof as that of Theorem \ref{thm:rateofconvgauss} using the fact that $\rho(\bar U)=1$ and Lemma \ref{lem:B100}, we have
\[
d_{Hell}(\rho^{\delta,0},\rho^{\delta,0,L})\le cL^{1/2}2^{-L}.
\]
From this and Theorem \ref{thm:rateofconvgauss} 
we get the conclusion. \eproof
\begin{remark}
We do not analyze approximation of the forward problem by choosing only a finite number of $J$ terms in the expansion of $K$. However, if the functions $\psi_j$ in \eqref{eq:uniform} and \eqref{eq:gaussian} have a decay rate $O(1/j^s)$ for $s>1$ with respect to the $L^\infty(D)$ norm, we can establish a rate of convergence in terms of $J$ for the approximated posterior measure which is similar to that considered in \cite{HoangScSt12}.
\end{remark}

%%%%%%%%%%%%%%%%%%%%%%%%%%%%%%%%%%%%%%
\section{Numerical examples}\label{sec:numerical}
\begin{figure}[h!]
	\centering
	\includegraphics[width=8cm]{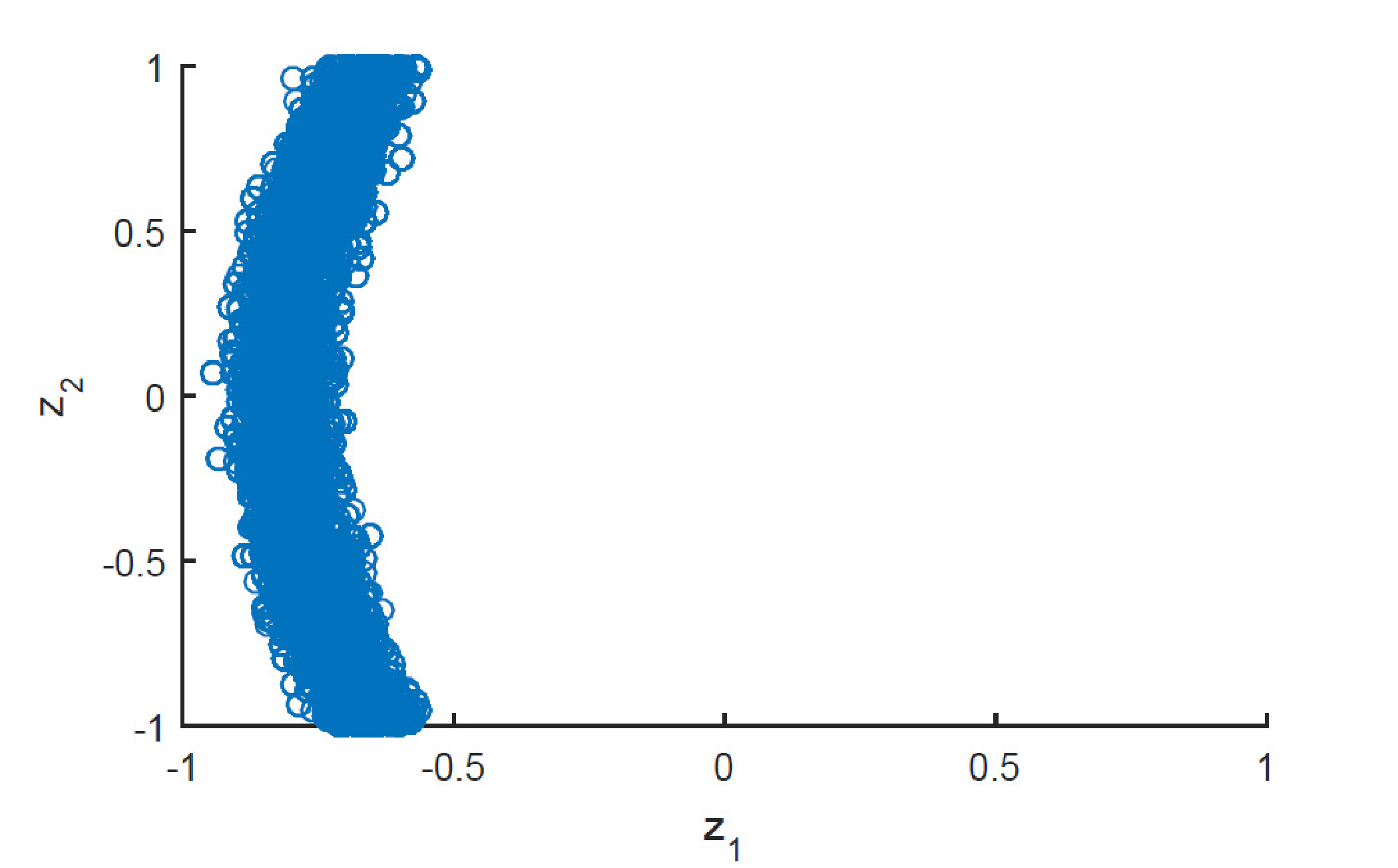}\\
	\captionof{figure}{60000 MCMC samples}
	\label{fig:1d_uni_u0_scatter}
\end{figure}

\begin{figure}[h!]
\centering
\begin{subfigure}{.5\textwidth}
  \centering
  \includegraphics[width=\textwidth]{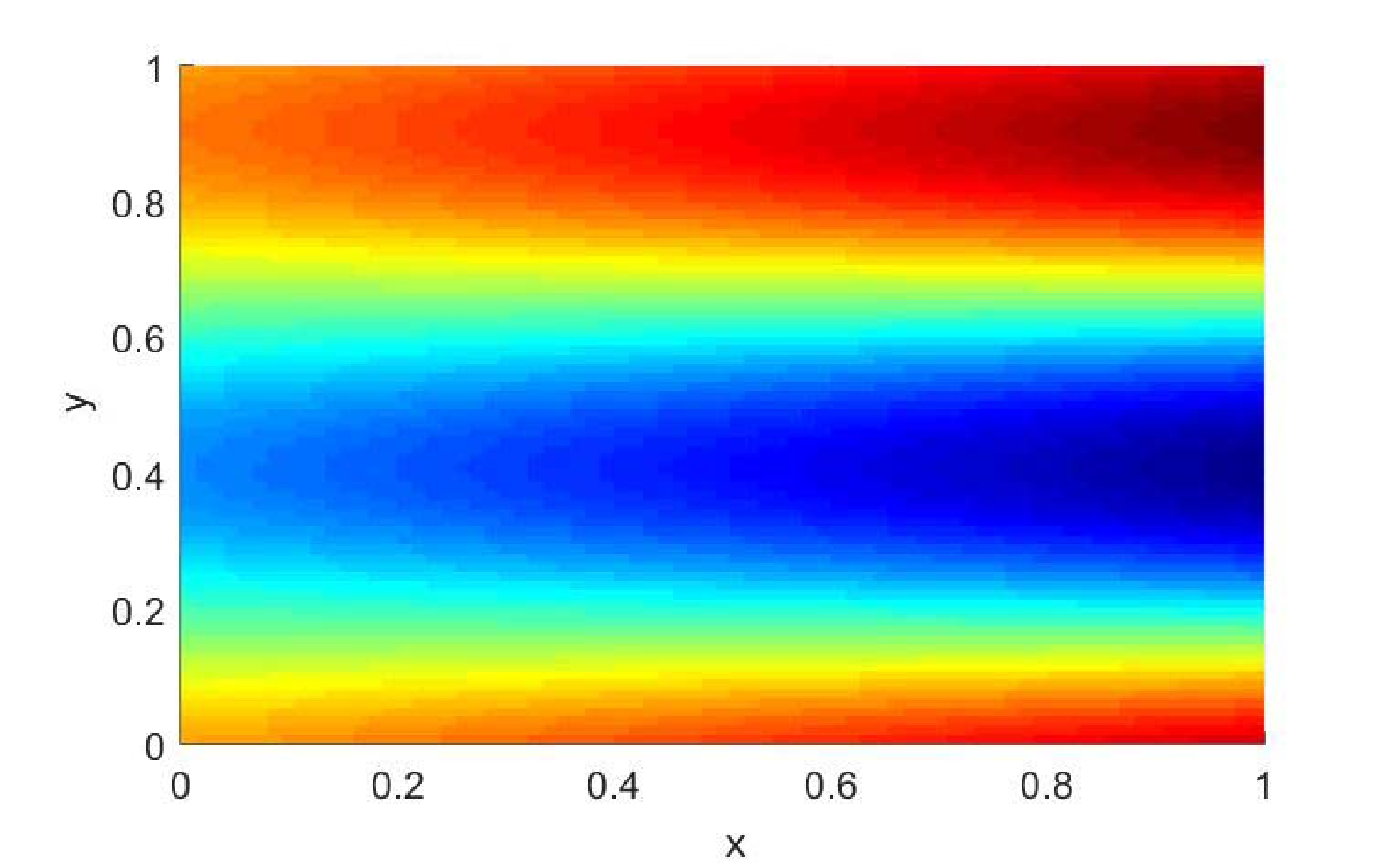}
  \captionof{figure}{Reference coefficient}
  \label{fig:1d_uni_u0_reference}
\end{subfigure}%
\begin{subfigure}{.5\textwidth}
  \centering
  \includegraphics[width=\textwidth]{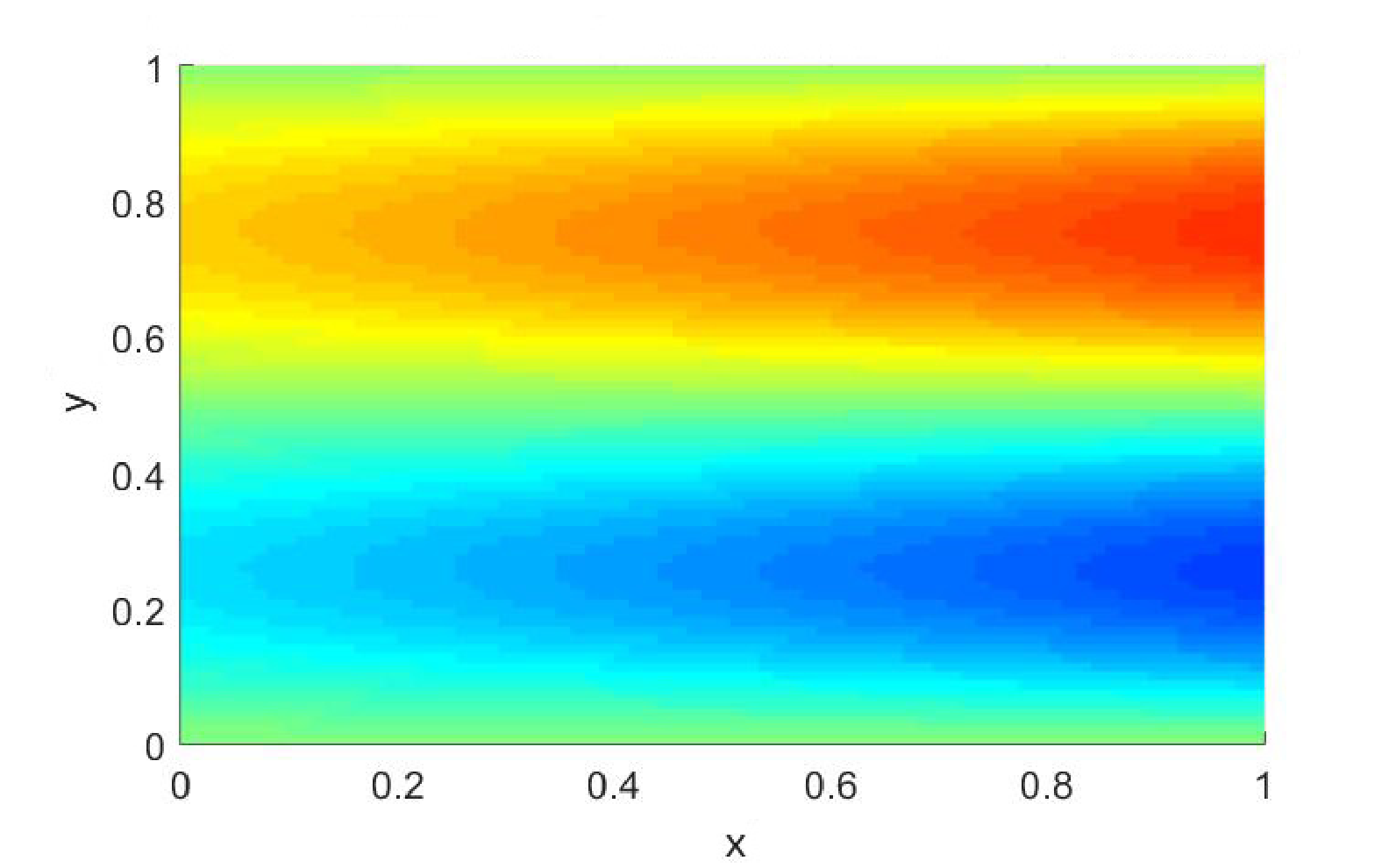}
  \captionof{figure}{Arithmetic average of the 60000 MCMC samples}
  \label{fig:1d_uni_u0_sample}
\end{subfigure}
\caption{The reference coefficient and the average of the MCMC samples}
\label{fig:1d_uni_u0}
\end{figure}

We perform the MCMC method in this section for some particular examples of recovering two scale coefficients. We use independent sampler for MCMC in all the numerical simulations in this section. We first consider the case of the  uniform prior probability measure. We consider a reference solution by taking a random realization of $z$ in \eqref{eq:uniform} or \eqref{eq:gaussian}. The data are obtained by adding a random realization of the noise into $\cO_i^0$ in \eqref{eq:Oi0}. We note that when $\ep$ is sufficiently small, this is a good approximation of the data generated in \eqref{eq:Oi}. 
%\subsection{1D Uniform with only observations for $u_0$}
%- Example for 1d with only observation for $u_0$, figure 1 shows the MCMC samples form a curve, figure 2 has two panels: left one show the reference, right one show the mean of mcmc sample\\
%
First for the one dimensional domain $D=(0,1)$, we consider the coefficient of the form
\be
A(z;x,y) = 9 + z_1(1+x)\sin(2\pi y) + z_2(1+x)\cos(2\pi y)
\label{eq:a1dunif}
\ee
where $z_1$ and $z_2$ are uniformly distributed in $[-1,1]$. We consider first the case where the functions $\ell_i$ in \eqref{eq:Oi} do not depend on $y$. As $\ue\wc u_0$ in $V$, we have that
\[
\lim_{\ep\to 0}\int_D\ell_i(x)\cdot\nabla\ue(x)dx=\int_D\ell_i(x)\cdot\nabla u_0(x)dx
\]
so that when $\ell_i$ do not depend on $y$, we essentially get information on $u_0$ only. 
Let $\ell_1(x,y)=x$ and $\ell_2(x,y)=x^2$. The observations become
\[
\cO_1^\ep=\int_Dx\nabla\ue(z;x)dx\approx\cO_1^0(z)=\int_D x\nabla u_0(z;x)dx,\ \ \cO_2^\ep=\int_Dx^2\nabla\ue(z;x)dx\approx\cO_2^0=\int_D x^2\nabla u_0(z;x)dx.
\]
The covariance of the noise is $10^{-3}I$ where $I$ is the $2\times 2$ identity matrix. We choose a reference pair $(z_1,z_2)$ at random. 
In Figure \ref{fig:1d_uni_u0_scatter}, we plot the first 60000 MCMC samples. We see that the samples are evenly scattered over a curve; the figure does not indicate clearly the value of the reference $(z_1,x_2)$. 
In Figure \ref{fig:1d_uni_u0}, we show the reference coefficient $A$ and the arithmetic average of 60000 MCMC samples. The figure shows that the average of the MCMC samples does not describe the reference coefficient accurately. 

\begin{figure}[h!]
	\centering
	\includegraphics[width=8cm]{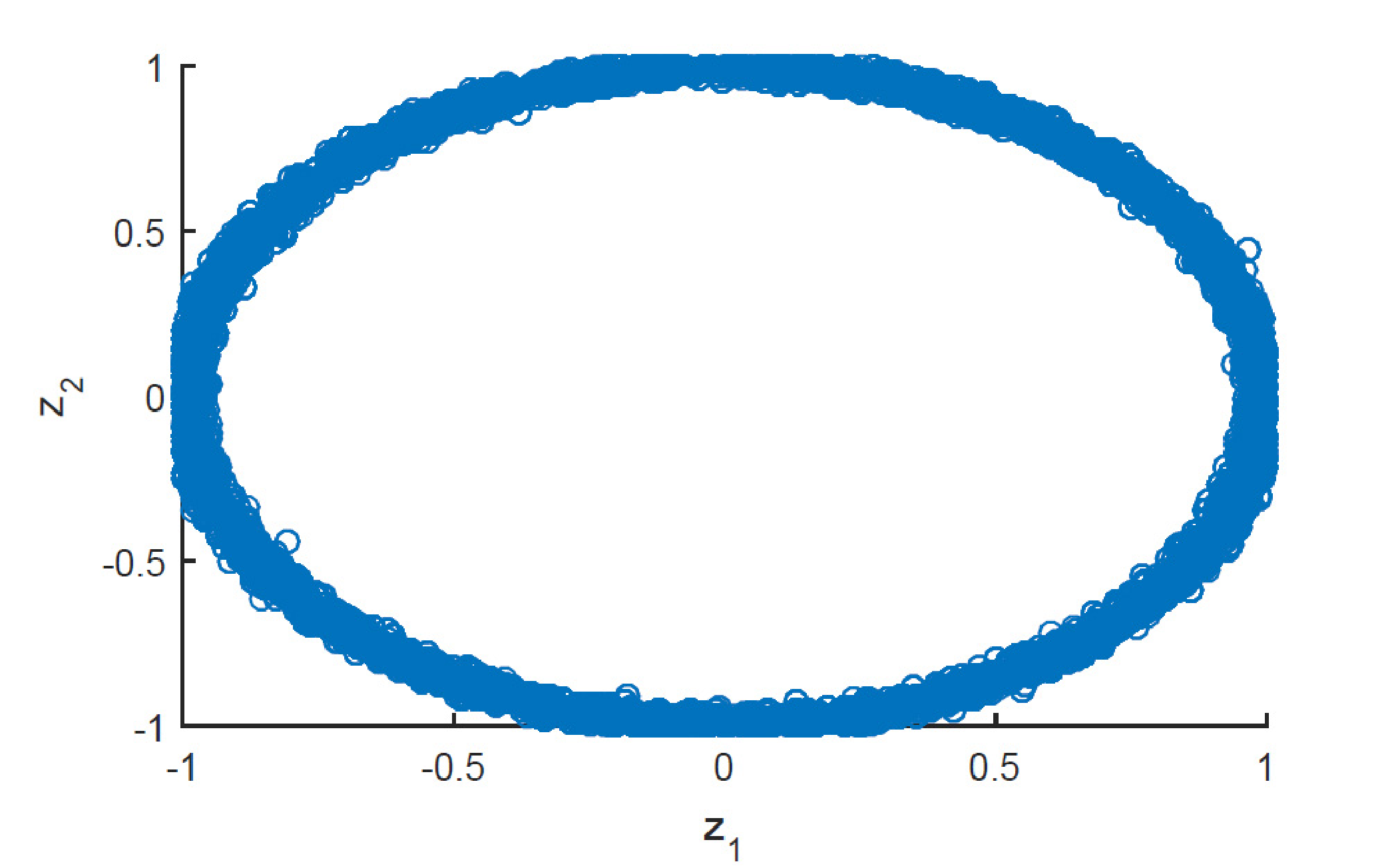}\\
	\captionof{figure}{60000 MCMC samples}
	\label{fig:1d_uni_lfa_scatter}
\end{figure}

\begin{figure}[h!]
\centering
\begin{subfigure}{.5\textwidth}
  \centering
  \includegraphics[width=\textwidth]{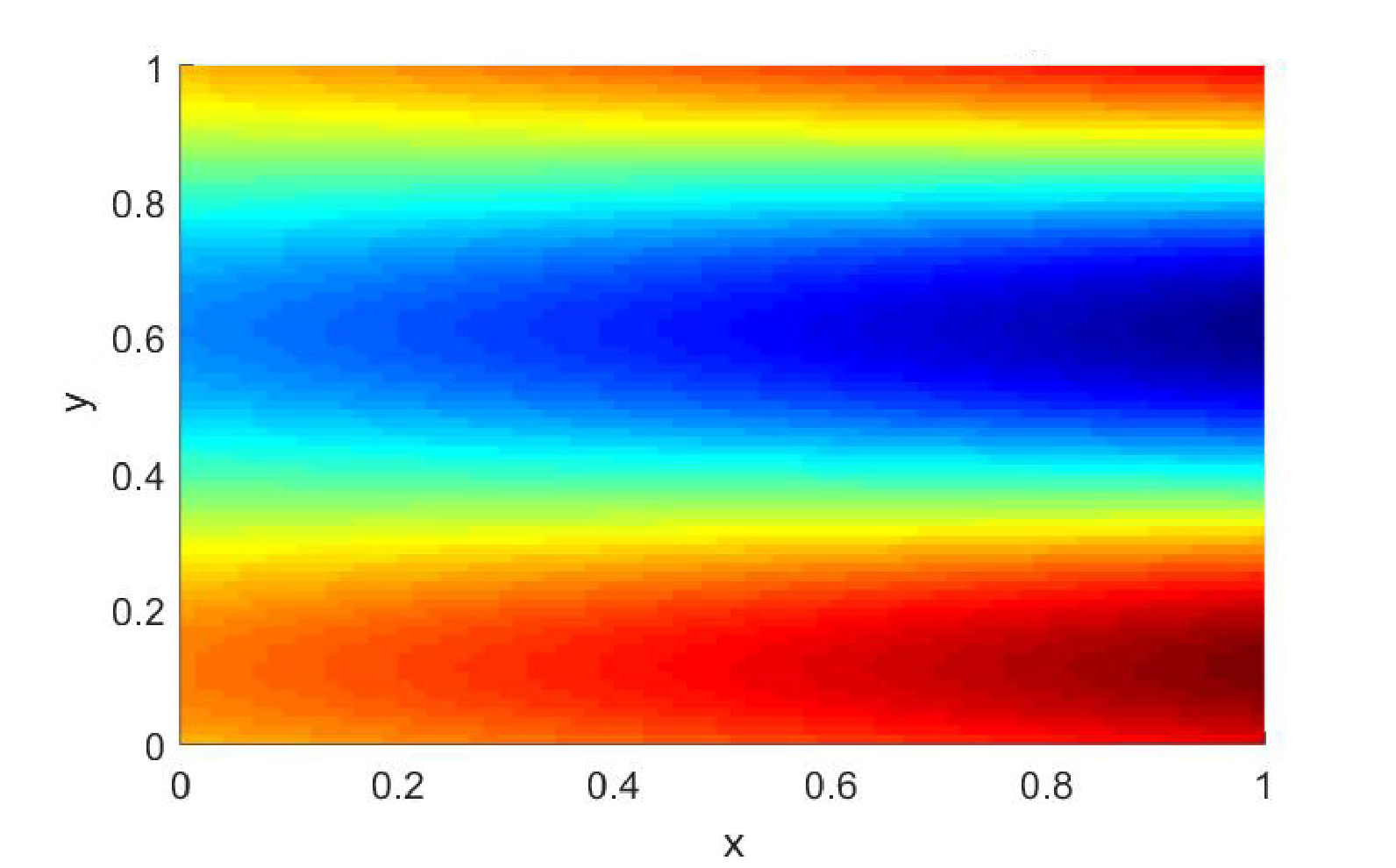}
  \captionof{figure}{Reference coefficient}
  \label{fig:1d_uni_lfa_reference}
\end{subfigure}%
\begin{subfigure}{.5\textwidth}
  \centering
  \includegraphics[width=\textwidth]{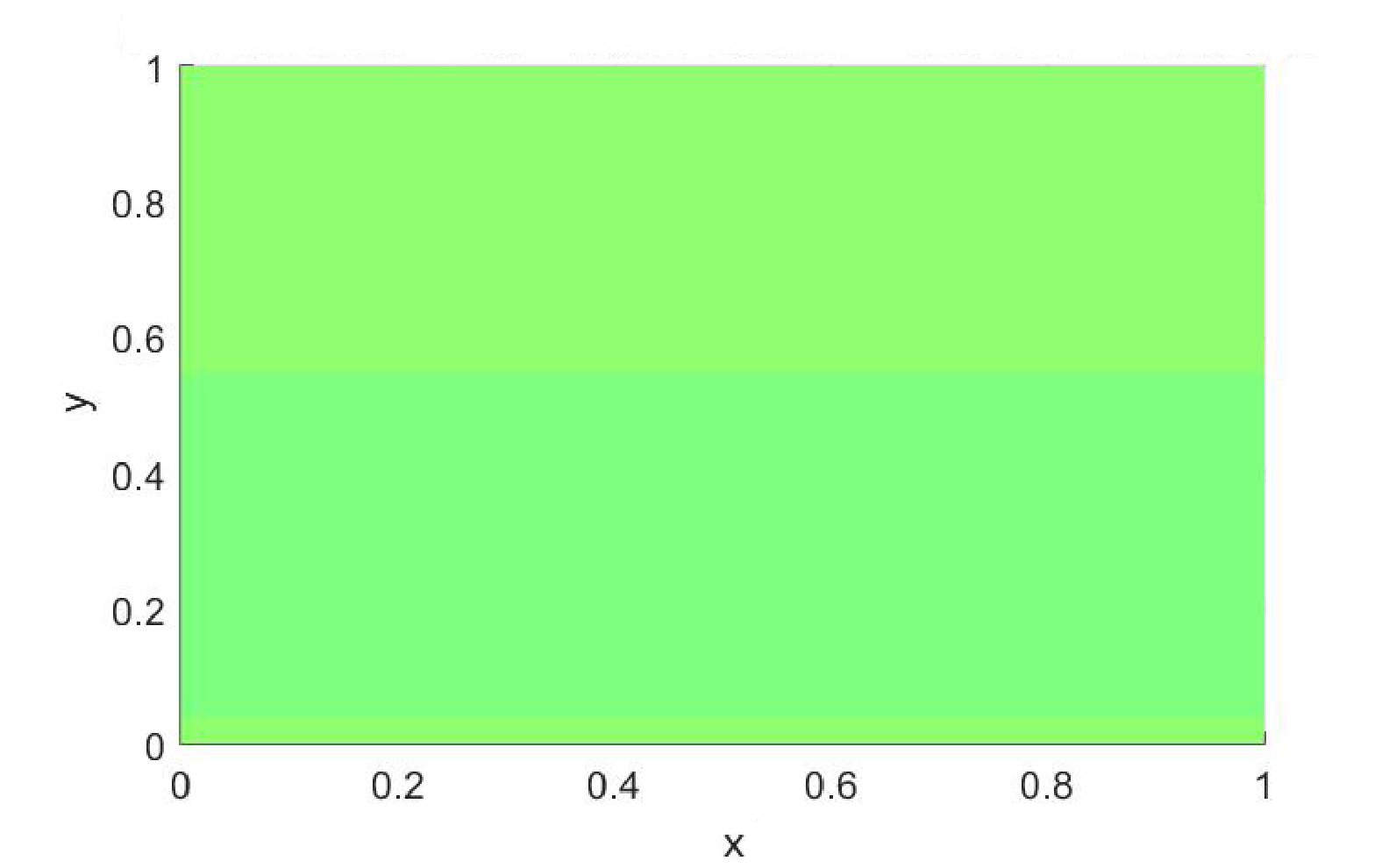}
  \captionof{figure}{Arithmetic average of the 60000 MCMC samples}
  \label{fig:1d_uni_lfa_sample}
\end{subfigure}
\caption{The reference coefficient and the average of the MCMC samples}
\label{fig:1d_uni_flux}
\end{figure}

Now we consider the case where the observation is  the flux of the two scale equation. This case can be analyzed in the same way as in  the previous sections. The coefficient is \eqref{eq:a1dunif}. 
%Coefficient : \\
%$$a(x,y,z) = 9 + z_1(1+x)\sin(2\pi y) + z_2(1+x)\cos(2\pi y)$$
%
We consider one observation
\[
\cO_1^\ep(z)=\int_DA(z;x,{x\over\ep})\nabla \ue(z;x)dx\approx\cO_1^0(z)=
\int_D\int_Y A(z;x,y)(\nabla u_0(z;x) + \nabla_y u_1(z;x,y))dydx
\]
Figure \ref{fig:1d_uni_lfa_scatter} shows the 60000 MCMC samples for $z=(z_1,z_2)$. The posterior does not provide good information on the reference $z=(z_1,z_2)$. Figure \ref{fig:1d_uni_flux} shows the reference coefficient and the arithmetic average of the coefficients obtained from the 60000 MCMC samples. We see that we cannot recover any details of the reference coefficient from the average  of the MCMC samples. We note that the flux of the two scale equation converges to the flux of the homogenized equation, i.e. 
\[
\lim_{\ep\to 0}\int_D A^\ep(z;x,{x\over\ep})\nabla\ue(z;x) dx=\int_D A^0(z;x)\nabla u_0(z;x)dx.
\]
Thus an observation on the flux essentially provides only the information on the homogenized equation, i.e. only the macroscopic information. This explains why we cannot recover accurately the reference coefficient. These examples show that it is necessary to have observations on $u_1$ to recover the microscopic structures. 

\begin{figure}[h!]
	\centering
	\includegraphics[width=8cm]{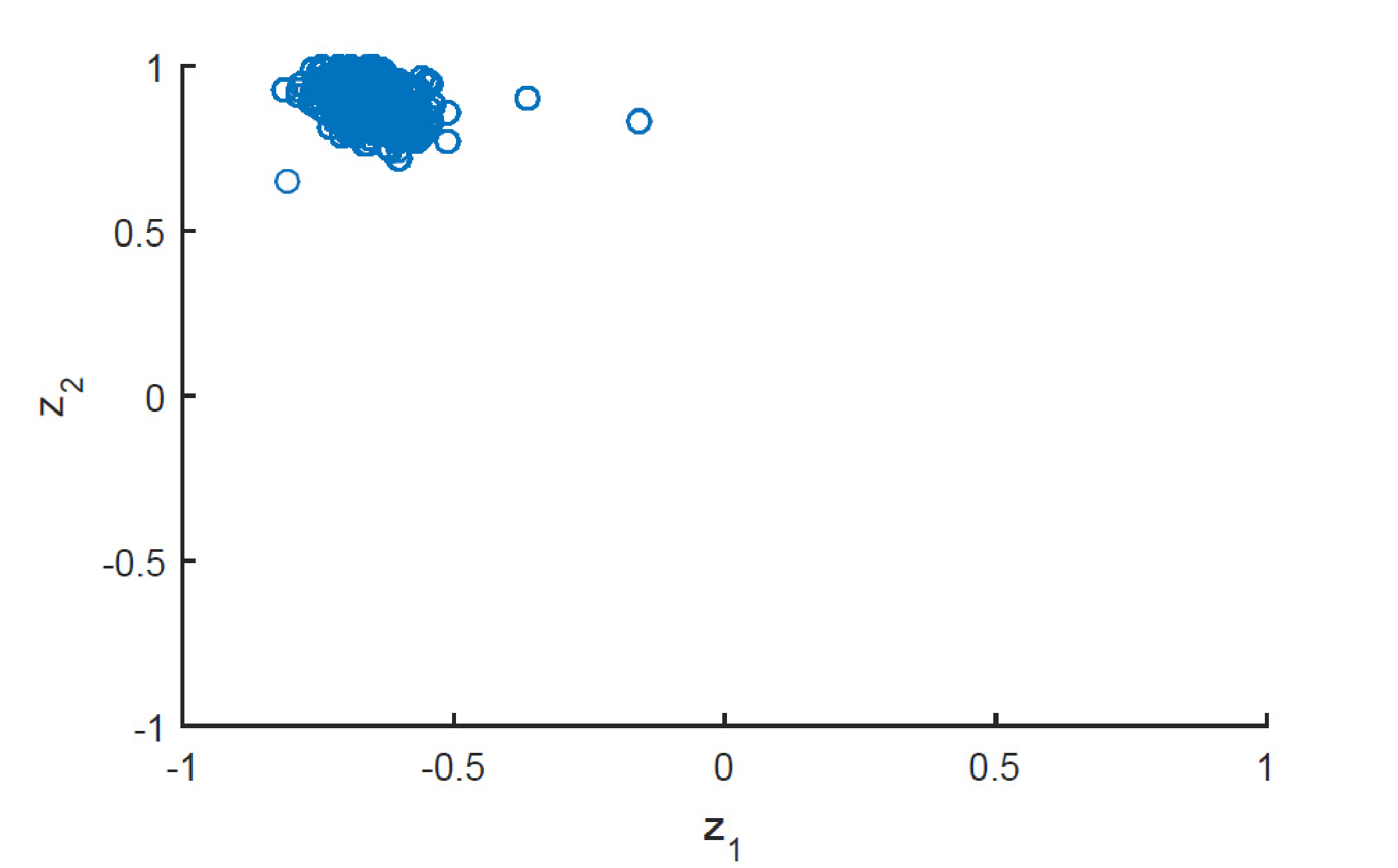}\\
	\captionof{figure}{60000 MCMC samples}
	\label{fig:1d_uni_u0u1_scatter}
\end{figure}

\begin{figure}[h!]
\centering
\begin{subfigure}{.5\textwidth}
  \centering
  \includegraphics[width=\textwidth]{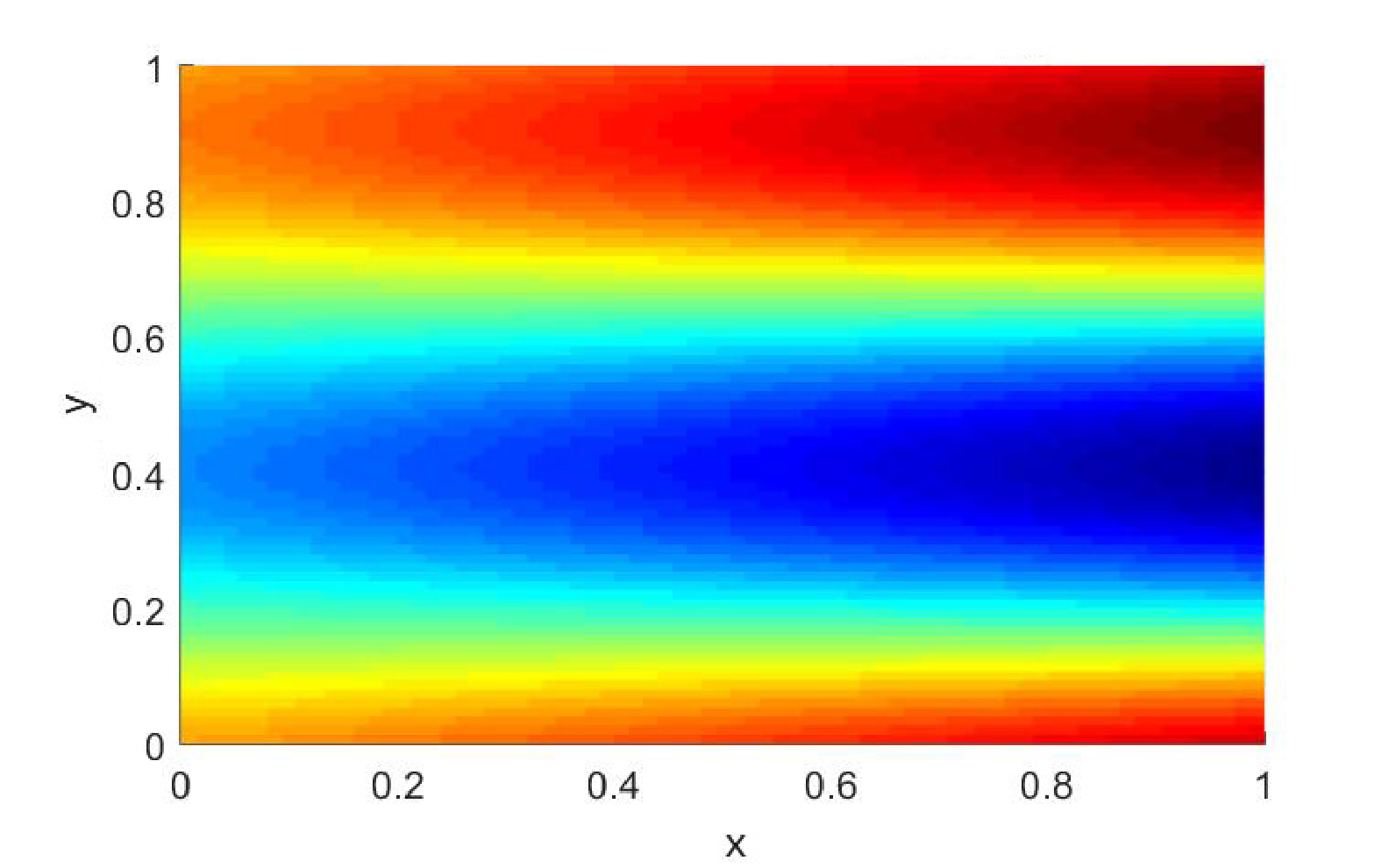}
  \captionof{figure}{Reference coefficient}
  \label{fig:1d_uni_u0u1_reference}
\end{subfigure}%
\begin{subfigure}{.5\textwidth}
  \centering
  \includegraphics[width=\textwidth]{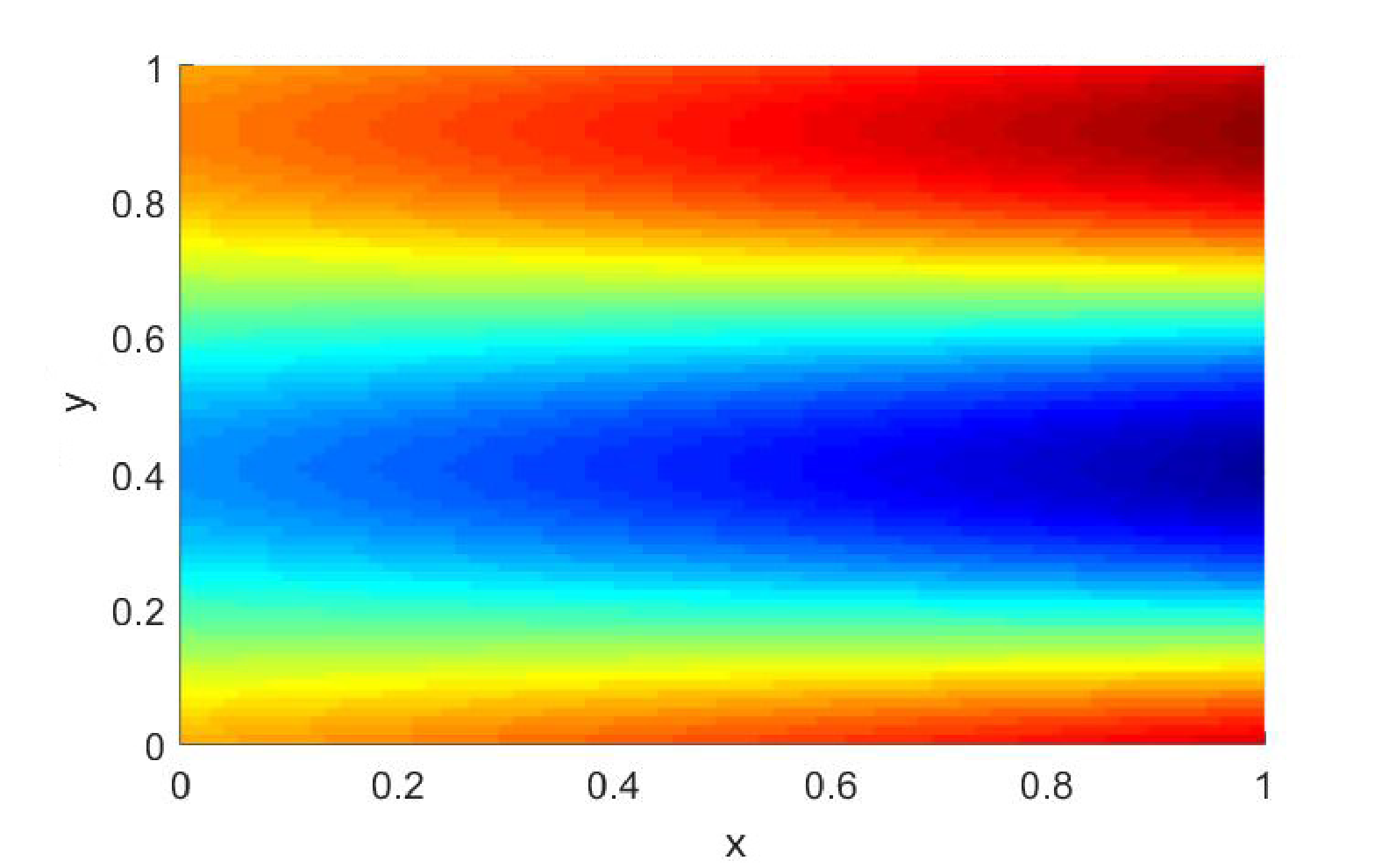}
  \captionof{figure}{Arithmetic average of the 60000 MCMC samples}
  \label{fig:1d_uni_u0u1_sample}
\end{subfigure}
\caption{The reference coefficient and the average of the MCMC samples}
\label{fig:1d_uni_u0u1}
\end{figure}

\begin{figure}[h!]
\centering
\begin{subfigure}{.5\textwidth}
  \centering
  \includegraphics[width=\textwidth]{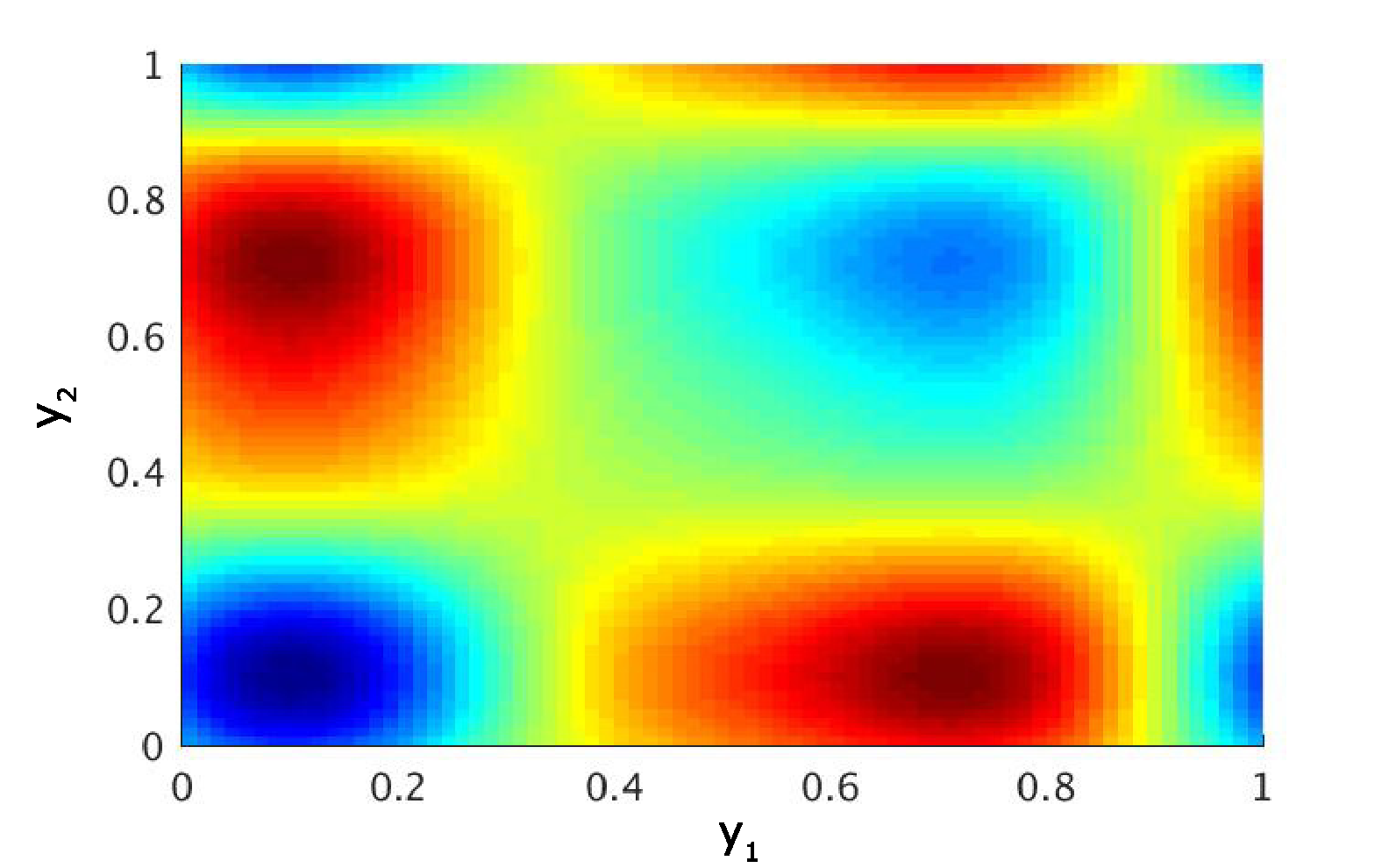}
  \captionof{figure}{Reference coefficient}
  \label{fig:2d_uni_reference}
\end{subfigure}%
\begin{subfigure}{.5\textwidth}
  \centering
  \includegraphics[width=\textwidth]{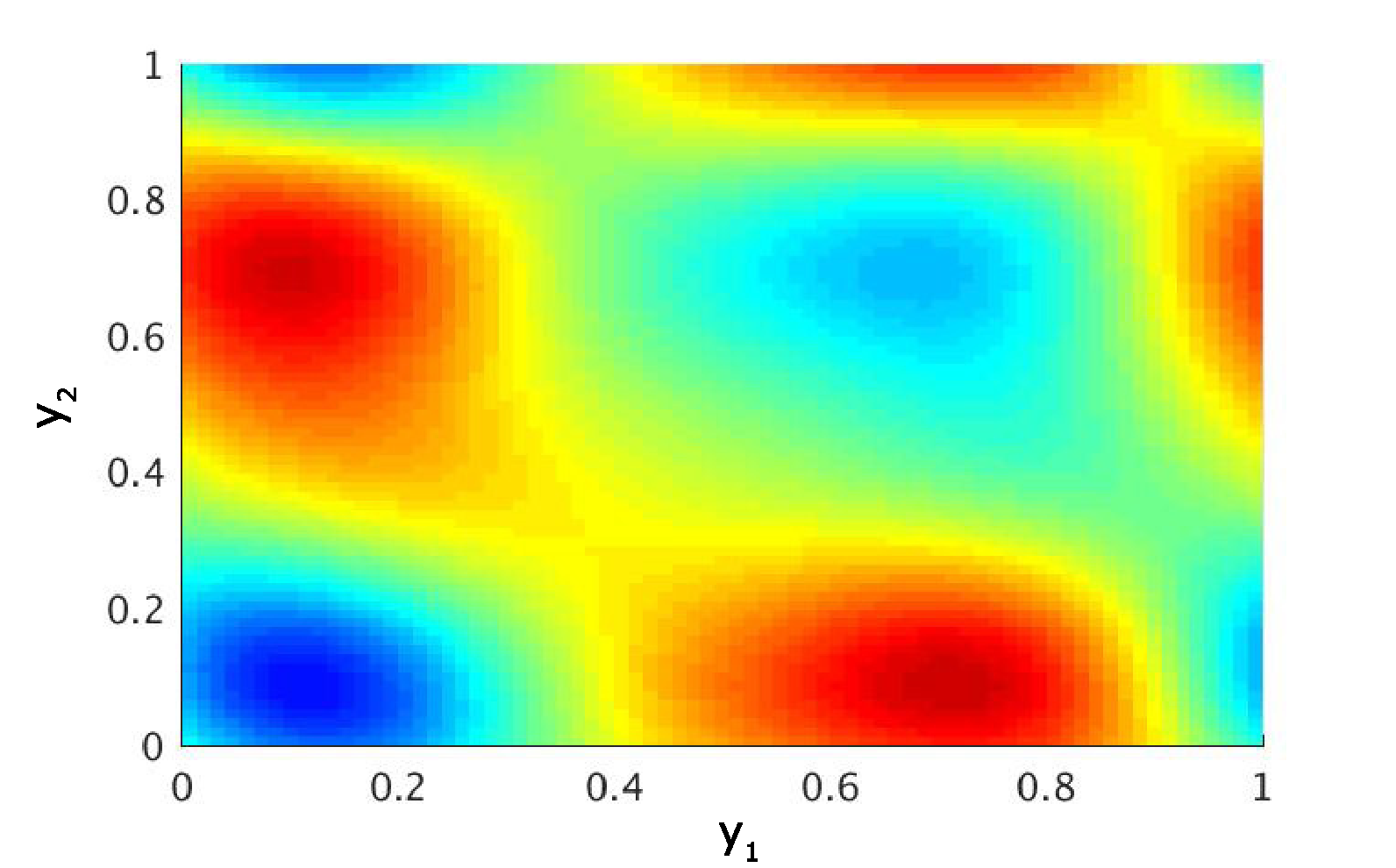}
  \captionof{figure}{Arithmetic average of the 120000 MCMC samples}
  \label{fig:2d_uni_sample}
\end{subfigure}
\caption{The reference coefficient and the average of the MCMC samples for $x=(0.2,0.2)$}
\label{fig:2d_uni}
\end{figure}

%\subsection{1D Uniform with observations for $u_0$ and $u_1$}
Next we consider the case where $\ell_i$ in \eqref{eq:Oi} depend also on $y$.
% Example for 1d with 2 observations for $u_0, u_1$, Figure 3 shows MCMC %samples, figure 4 shows reference and mean of mcmc samples\\
%
%Coefficient : \\
%$$a(x,y,z) = 9 + z_1(1+x)\sin(2\pi y) + z_2(1+x)\cos(2\pi y)$$
The coefficient $A$ is of the form \eqref{eq:a1dunif}. The functions
\[
\ell_1(x,y)=x(1+\sin(2\pi y)),\ \ \ell_2(x,y)=x(1+\cos(2\pi y)).
\]
The observations are:
\beqas
\cO_1^\ep(z)=\int_Dx(1+\sin(2\pi{x\over\ep}))\nabla\ue(x)dx\approx\cO_1^0(z)=\int_D\int_Y [x(1+\sin(2\pi y))](\nabla u_0 + \nabla_y u_1)dydx,\\
\cO_2^\ep(z)=\int_Dx(1+\cos(2\pi{x\over\ep})\nabla\ue(x)dx\approx\cO_2^0(z)=\int_D\int_Y [x(1+\cos(2\pi y))](\nabla u_0 + \nabla_y u_1)dydx.
\eeqas
In Figure \ref{fig:1d_uni_u0u1_scatter} we plot the 60000 MCMC samples for the pair $(z_1,z_2)$. The figure shows that the posterior measure provides very good prediction on the reference $z=(z_1,z_2)$. In Figure \ref{fig:1d_uni_u0u1}, we show the reference coefficient and the average of the coefficient obtained from the MCMC samples. The figure shows  that we have a good recovery of the reference coefficient.  
\\
Now we consider a two dimensional problem in the domain $D=(0,1)^2$. Here $x=(x_1,x_2)\in D$ and $y=(y_1,y_2)\in Y=(0,1)^2$ -- the unit cube in $\IR^2$. We consider the coefficient of the form
\[
a(x,y,z) = 10 + \frac{1}{4}(1+x_1)(1+x_2)\sum z_{F_1,F_2}F_1(y_1)F_2(y_2)
\]
where the summation is over all %$i,j=1,2$;
$F_i(y_i)\in \{\sin(2\pi y_i), \cos(2\pi y_i), \frac14\sin(4\pi y_i), \frac14\cos(4\pi y_i)\}$, $i=1,2$. The random variables $z_{F_1,F_2}$ are uniformly distributed in $[-1,1]$. We consider the observations of the form
\[
\int_D(1+x_1)(1+x_2)G_1({x_1\over\ep})G_2({x_2\over\ep}){\partial\ue\over\partial x_p}dx\approx\int_D\int_Y (1+x_1)(1+x_2)G_1(y_1)G_2(y_2)({\partial u_0\over\partial x_p}(z;x) + {\partial u_1\over\partial y_p}(z;x,y))dy dx
\]
for all $G_i(y_i)\in \{\sin(2\pi y_i), \cos(2\pi y_i), \sin(4\pi y_i), \cos(4\pi y_i), \sin(6\pi y_i), \cos(6\pi y_i)\}$, $i=1,2$, $p=1,2$. Figure \ref{fig:2d_uni} presents the reference coefficient and the average of 120000 MCMC samples fixing $x=(0.2,0.2)$. It shows that we have a reasonably good recovery of the two scale coefficient. 
\\

\begin{figure}[h!]
\centering
\begin{subfigure}{.5\textwidth}
  \centering
  \includegraphics[width=\textwidth]{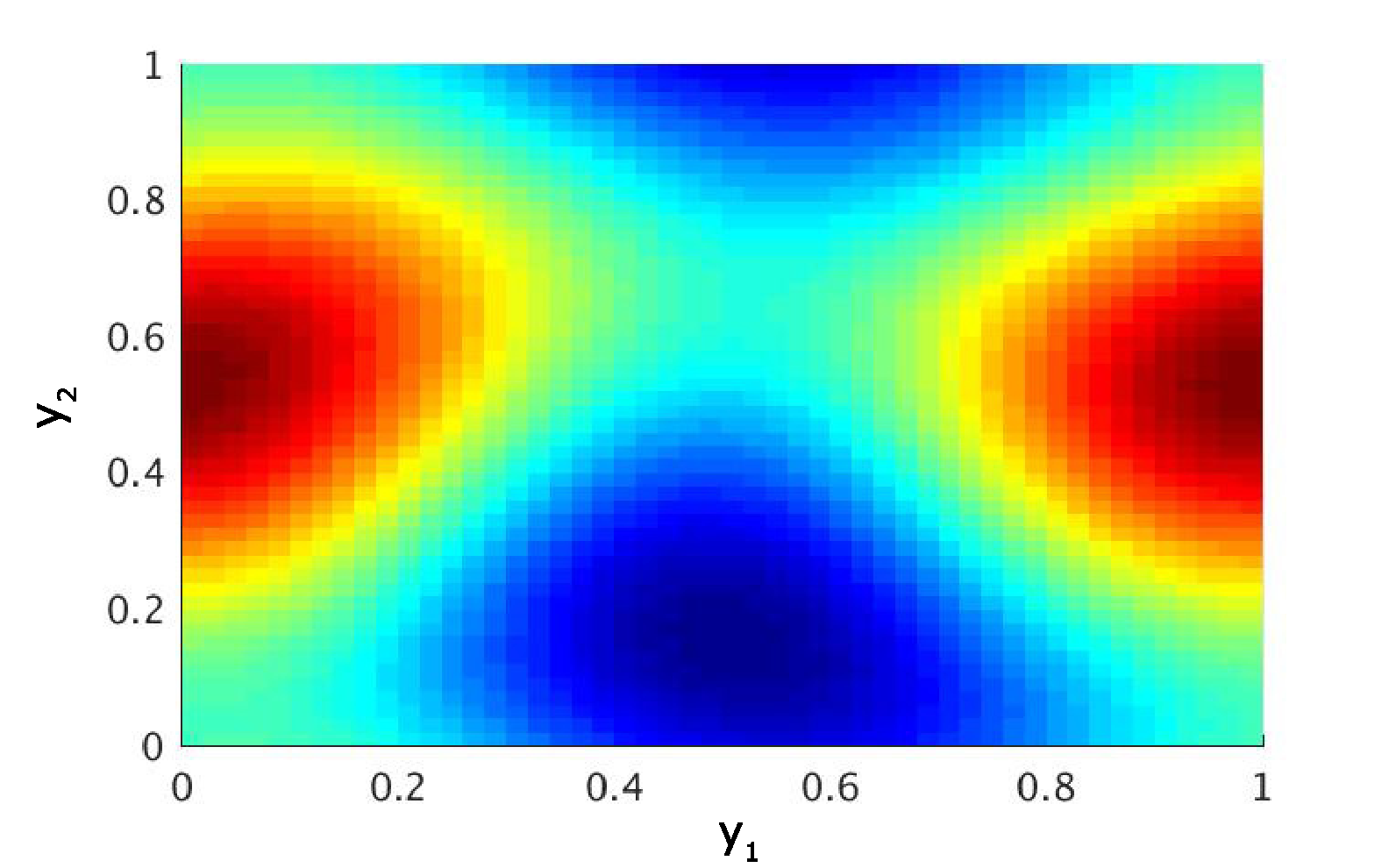}
  \captionof{figure}{Reference coefficient}
  \label{fig:2d_logn_reference2525}
\end{subfigure}%
\begin{subfigure}{.5\textwidth}
  \centering
  \includegraphics[width=\textwidth]{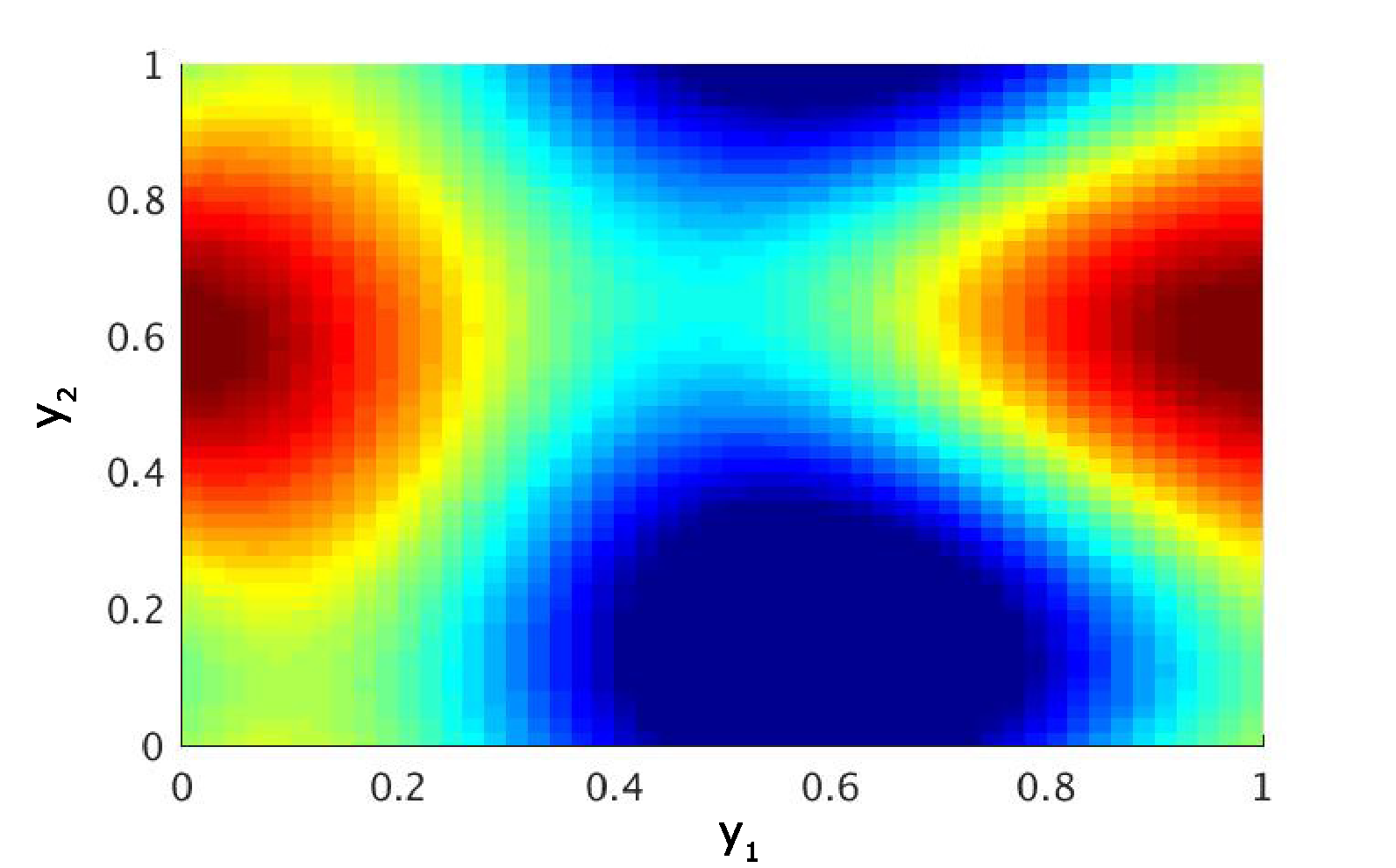}
  \captionof{figure}{Arithmetic average of the 120000 MCMC samples}
  \label{fig:2d_logn_sample2525}
\end{subfigure}
\caption{Reference coefficient and the average  of MCMC samples for $x=(0.25,0.25)$}
\label{fig:2d_logn_80z_2525}
\end{figure}

\begin{figure}[h!]
\centering
\begin{subfigure}{.5\textwidth}
  \centering
  \includegraphics[width=\textwidth]{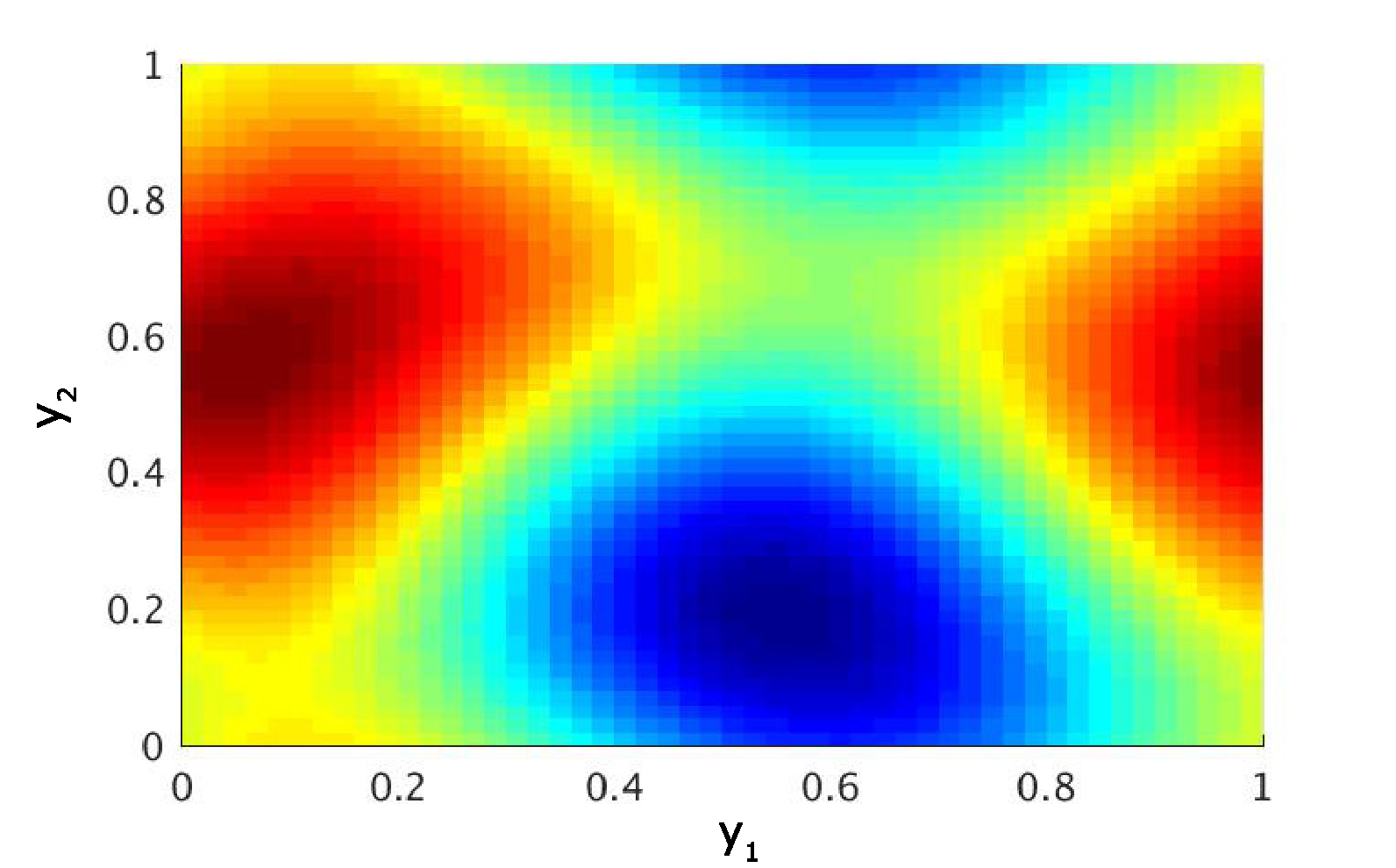}
  \captionof{figure}{Reference coefficient}
  \label{fig:2d_logn_reference2575}
\end{subfigure}%
\begin{subfigure}{.5\textwidth}
  \centering
  \includegraphics[width=\textwidth]{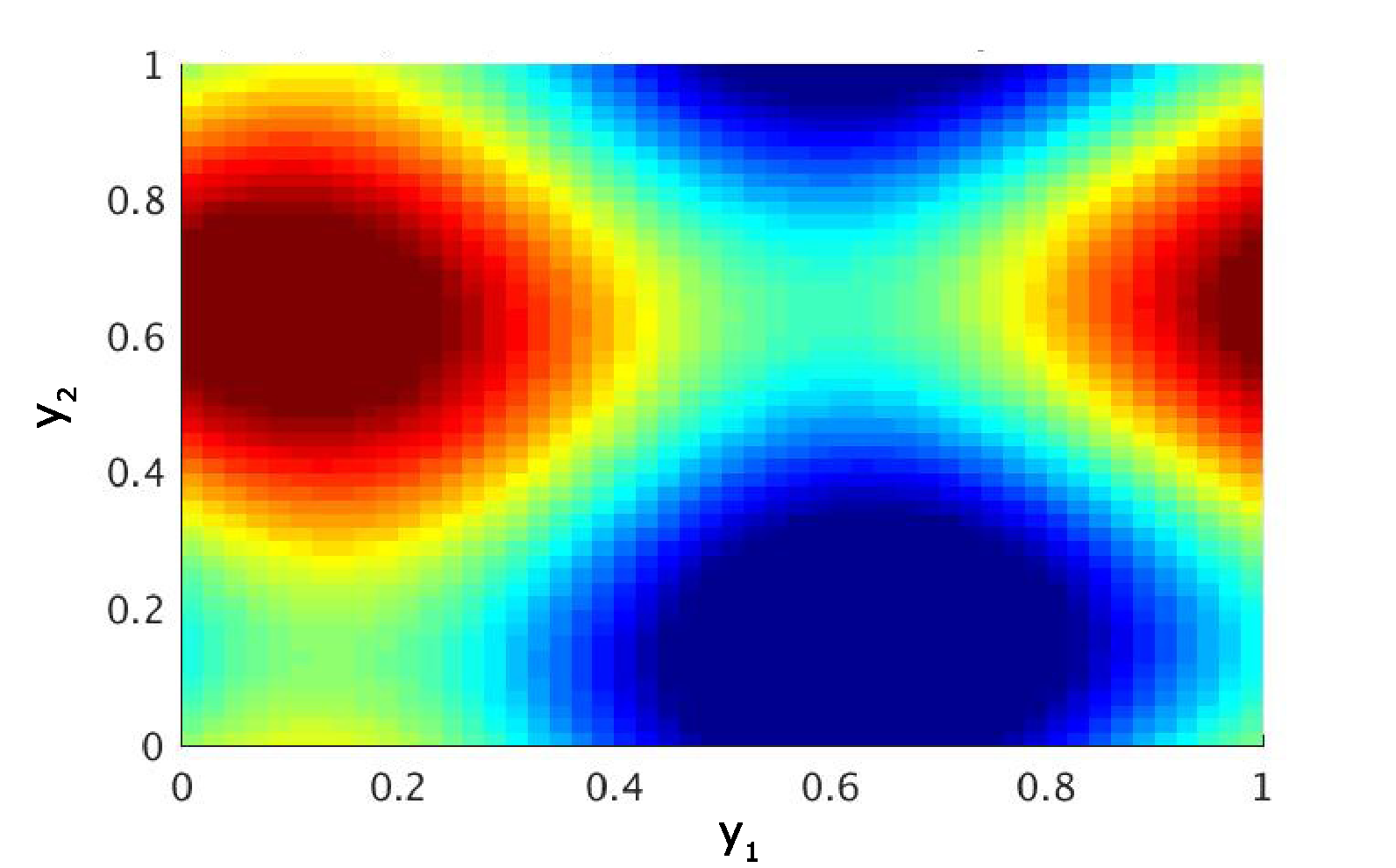}
  \captionof{figure}{Arithmetic average of the 120000 MCMC samples}
  \label{fig:2d_logn_sample2575}
\end{subfigure}
\caption{Reference coefficient and the average  of MCMC samples for $x=(0.25,0.75)$}
\label{fig:2d_logn_80z_2575}
\end{figure}

%\subsection{2D Lognormal}
%- Example for 2d logn\\
We now consider the case of an equation in the two dimensional domain $D=(0,1)^2$ with the log-Gaussian prior. 
%Coefficient :\\
%\[
%a(x,y,z) = \exp\left(  \sum_{k_1+k_2+l_1+l_2 \neq 0 }z %\frac{H_a(x_1;k_1)H_b(x_2;k_2)H_c(y_1;l_1)H_d(y_2;l_2)}{( (2\pi k_1)^2 + %(2\pi k_2)^2 + (2\pi l_1)^2 + (2\pi l_2)^2 )^2}  \right)
%\]
%where $H_i(t;k) \in \{ \cos(2\pi kt) | k=0 \} \bigcup \{ \sin(2\pi kt) | %k=1 \} \bigcup \{ \cos(2\pi kt) | k=1 \}$\\
Let $\psi_i(x)$, $i=1,\ldots,9$ be the 9 eigenfunctions of $-\Delta_x^{-1}$ with periodic boundary condition on $D=(0,1)^2$ that are of the form $\sin(\cos)(2k\pi x_1)\sin(\cos)(2l\pi x_2)$ with $k,l=0,1$ (excluding the zero function), with eigenvalues $\lambda_i$. Let $\phi_j(y)$, $j=1,\ldots,9$, be the 9 eigenfunctions of $-\Delta_y^{-1}$ with periodic boundary condition on $Y=(0,1)^2$ that are of the form $\phi_j(x)=\sin(\cos)(2k\pi y_1)\sin(\cos)(2l\pi y_2)$ with $k,l=0,1$ (excluding the zero function) with eigenvalues $\mu_j$. 
%$\psi_i(x)\in H_1(x_1)H_2(x_2)$ where $H_a(x_a)\in \{1, \cos(2\pi x_a), \sin(2\pi x_a)\}$, and similarly for $\phi_j$. 
We let $\psi_1(x)=1$, $\lambda_1=0$ and $\phi_1(y)=1$ with $\mu_1=0$. 
We consider the coefficient
\[
a(z;x,y)=\exp\left(\sum_{i,j=1,\ldots,9\atop (i,j)\ne (1,1)}z_{i,j}{1\over(\lambda_i+\mu_j)^2}\psi_i(x)\psi_j(y)\right)
\]
where $z_{ij}\sim {\cal N}(0,1)$. 
All together we have 80 terms in the summation. We consider all the observations of the form
\beqas
&&\int_D1000F_1(x_1)F_2(x_2)G_1({x_1\over\ep})G_2({x_2\over\ep}){\partial\ue\over\partial x_p}(z;x)dx\approx\\
&&\int_D\int_Y 1000F_1(x_1)F_2(x_2)G_1(y_1)G_2(y_2)( {\partial u_0\over\partial x_p}(z;x) + {\partial u_1\over\partial y_p}(z;x,y))dydx 
\eeqas
for all
$ F_i(x_i) \in \{1, 1 + \sin(2\pi kx_i), 1 + \cos(2\pi k x_i), k=1,2 \} $ and $G_i(y_i) \in \{1, 1 + \sin(2\pi ky_i), 1 + \cos(2\pi k y_i), k=1,2\}$; $p=1,2 $.
In Figure \ref{fig:2d_logn_80z_2525} we show the reference coefficient for $x=(0.25,0.25)$ and the average of 120000 MCMC samples. The figure shows that the MCMC method provides a reasonably good recovery of the coefficient.  
Simimarly, for $x=(0.25, 0.75)$ we show the reference coefficient and the average of the MCMC samples in Figure \ref{fig:2d_logn_80z_2575}. 
\\
\\
{\bf Acknowledgement} The research is supported by the Singapore MOE AcRF Tier 1 grant RG30/16, the MOE Tier 2 grant MOE2017-T2-2-144, and a graduate scholarship from Nanyang Technologial University, Singapore. 
\begin{appendices}
\section{}\label{app:a}

We show Theorem \ref{thm:rateofconvunif} in this appendix. We show that $I_1$ and $I_2$ in \eqref{eq:I1} and \eqref{eq:I2} have upper bound $c\ep$. We first recall the following result.
\begin{lemma}\label{lem:1}
For the coefficient of the form \eqref{eq:uniform}, under Assumption \ref{assum:furtherregularitypsiunif}, there is a constant $c$ independent of $z\in U$ such that
\[
\|\nabla\ue(z)-[\nabla u_0(z)+\nabla_yu_1(z;\cdot,{\cdot\over\ep})\|_{L^2(D)}\le c\ep^{1/2}.
\]
\end{lemma}
\bproof
This result is indeed Theorem 5.2 of Hoang and Schwab \cite{HSmultirandom}. From Assumption \ref{assum:furtherregularitypsiunif}, we have that $A(z;x,y)$ is uniformly bounded in $C^1(\bar D,C^{1,1}(\bar Y))$. The result follows from Theorem 5.2 of \cite{HSmultirandom}.\eproof

From Lemma \ref{lem:1}, we have that
\be
\left|\int_D\ell_i(x,{x\over\ep})\cdot\nabla\ue(z) dx-\int_D\ell_i(x,{x\over\ep})\cdot(\nabla u_0(z;x)+\nabla_yu_1(z;x,{x\over\ep}))dx\right|\le c\ep^{1/2}.
\label{eq:a1}
\ee
Next we show the following result.
\begin{lemma}\label{lem:2}
Under Assumption \ref{assum:furtherregularitypsiunif}, when $\ell_i\in C^{0,1}(\bar D,C(\bar Y))^d$, we have
\be
\left|\int_D\ell_i(x,{x\over\ep})\cdot(\nabla u_0(z;x)+\nabla_yu_1(z;x,{x\over\ep}))dx-\int_D\int_Y\ell_i(x,y)\cdot(\nabla u_0(z;x)+\nabla_yu_1(z;x,y))dydx\right|\le c\ep^{1/2}.
\label{eq:a2}
\ee
\end{lemma}
\bproof
From Proposition 4.2 of \cite{HSmultirandom}, we have that $u_0(z;\cdot)$ is uniformly bounded in $H^2(D)$ with respect to $z$.
Let $D_m$, $m=1,\ldots,M$ be the cubes of the form
%of size $\ep$ which is obtained from $Y$ by translating by a vector whose coordinates are all intergers and rescaling by $\ep$, 
$\ep{\boldsymbol n}+\ep Y$ where ${\boldsymbol n}\in \IN^d$,
that are entirely contained in $D$. We denote by $D^\ep=\bigcup_{m=1}^M D_m$. We have that
\be
\left|\int_{D\setminus D^\ep}\ell_i(x,{x\over\ep})\cdot \nabla u_0(z;x)dx\right|\le c\int_{D\setminus D^\ep}|\nabla u_0(z;x)|dx\le c|D\setminus D^\ep|^{1/2}\left(\int_{D\setminus D^\ep}|\nabla u_0(z;x)|^2dx\right)^{1/2}\le c\ep^{1/2}. 
\label{eq:a3}
\ee
 From the proof of Lemma 5.5 in \cite{HSmultirandom}, we have
\be
\int_D\int_Y\left|{\partial u_0\over\partial x_k}(z;x)-{\partial u_0\over\partial x_k}\left(z;\ep\left[{x\over\ep}\right]+\ep t\right)\right|dtdx\le c\ep,
\label{eq:a4}
\ee
where $c$ is independent of $z$; ${\partial u_0\over\partial x_k}$ is continuously extended outside $D$. We therefore have
\be
\left|\int_{D^\ep}\ell_i(x,{x\over\ep})\cdot\nabla u_0(z;x)-\int_{D^\ep}\ell_i(x,{x\over\ep})\cdot\left(\int_Y\nabla u_0\left(z;\ep\left[{x\over\ep}\right]+\ep t\right)dt\right) dx\right|\le c\ep.
\label{eq:a5}
\ee
%Let $D_m$ $m=1,\ldots,M$ be the cube of size $\ep$ which is obtained from $Y$ by translating by a vector whose coordinates are all intergers and rescaling by $\ep$, that are entirely contained in $D$. We denote by $D^\ep=\bigcup_{m=1}^M D_m$. We have that
%\[
%\left|\int_{D\setminus D^\ep}\ell_i(x,{x\over\ep})\cdot \nabla u_0(x)dx\right|\le c\int_{D\setminus D^\ep}|\nabla u_0(x)|dx\le c|D\setminus D^\ep|\left\int_{D\setminus D^\ep}|\nabla u_0(x)|^2dx\right)^{1/2}\le c\ep^{1/2}. 
 We have
\beqas
 \int_{D^\ep}\ell_i(x,{x\over\ep})\cdot\left(\int_Y\nabla u_0\left(z;\ep\left[{x\over\ep}\right]+\ep t\right)dt\right) dx=
 \sum_{m=1}^M\int_{D_m}\ell_i(x,{x\over\ep})\cdot\left(\int_Y\nabla u_0\left(z;\ep\left[{x\over\ep}\right]+\ep t\right)dt\right) dx.
 \eeqas
Let $x_m$ be the centre of $D_m$. We have 
\begin{eqnarray}
&&\Big|\sum_{m=1}^M\int_{D_m}\ell_i(x,{x\over\ep})\cdot\left(\int_Y\nabla u_0\left(z;\ep\left[{x\over\ep}\right]+\ep t\right)dt\right) dx-
\sum_{m=1}^M\int_{D_m}\ell_i(x_m,{x\over\ep})\cdot\left(\int_Y\nabla u_0\left(z;\ep\left[{x\over\ep}\right]+\ep t\right)dt\right) dx\Big|\nonumber\\
&&\le \sum_{m=1}^Mc\ep\int_{D_m}\left|\int_Y\nabla u_0\left(z;\ep\left[{x\over\ep}\right]+\ep t\right)dt\right|\nonumber\\
&& =c\ep\int_{D^\ep}\left|\int_Y\nabla u_0\left(z;\ep\left[{x\over\ep}\right]+\ep t\right)dt\right| dx\le c\ep.
\label{eq:a6}
\end{eqnarray}
As $\int_Y\nabla u_0\left(z;\ep\left[{x\over\ep}\right]+\ep t\right)dt$ is constant for $x\in D_m$,
with a simple change of variable, for $x\in D_m$
\begin{eqnarray*}
\int_{D_m}\ell_i(x_m,{x\over\ep})\cdot \left(\int_Y\nabla u_0\left(z;\ep\left[{x\over\ep}\right]+\ep t\right)dt\right) dx=\ep^d\int_Y\nabla u_0\left(z;\ep\left[{x\over\ep}\right]+\ep t\right)dt\cdot\int_Y\ell_i(x_m,y)dy=\\
\int_{D_m}\int_Y\ell_i(x_m,y)\cdot\left(\int_Y\nabla u_0\left(z;\ep\left[{x\over\ep}\right]+\ep t\right)dt\right)dy dx
\eeqas
As for $x\in D_m$, $|\ell_i(x,y)-\ell_i(x_m,y)|\le c\ep$, from this equality, we have 
\beqas
&&\left|\int_{D_m}\ell_i(x_m,{x\over\ep})\cdot \left(\int_Y\nabla u_0\left(z;\ep\left[{x\over\ep}\right]+\ep t\right)dt\right) dx-\int_{D_m}\int_Y\ell_i(x,y)\cdot\left(\int_Y\nabla u_0\left(z;\ep\left[{x\over\ep}\right]+\ep t\right)dt\right) dydx\right|\\
&&\le c\ep\ep^d\int_Y|\nabla u_0\left(z;\ep\left[{x\over\ep}\right]+\ep t\right)|dt\\
&&=c\ep\int_{D_m}|\nabla u_0(z;x)|dx
\eeqas
as $\ep[x/\ep]+\ep t\in D_m$ for $x\in D_m$ and $t\in Y$, the last equality is obtained from a simple change of variable. Thus 
\begin{eqnarray}
&&\left|\sum_{m=1}^M\int_{D_m}\ell_i(x_m,{x\over\ep})\cdot \left(\int_Y\nabla u_0\left(z;\ep\left[{x\over\ep}\right]+\ep t\right)dt\right) dx-\int_{D^\ep}\int_Y \ell_i(x,y)\cdot\left(\int_Y\nabla u_0\left(z;\ep\left[{x\over\ep}\right]+\ep t\right)dt\right)dy dx\right|\nonumber\\
&&\le c\ep\int_{D^\ep}|\nabla u_0(z;x)|dx\le c\ep.
 \label{eq:a7}
\end{eqnarray}
From \eqref{eq:a4}, we have
\begin{eqnarray}
&&\left|\int_{D^\ep}\int_Y \ell_i(x,y)\cdot\left(\int_Y\nabla u_0\left(z;\ep\left[{x\over\ep}\right]+\ep t\right)dt\right)dy dx-\int_{D^\ep}\int_Y \ell_i(x,y)\cdot\nabla u_0(z;x) dydx\right|\nonumber\\
&&=\left|\int_{D^\ep}\left(\int_Y\ell_i(x,y)dy\right)\cdot\int_Y\left(\nabla u_0\left(z;\ep\left[{x\over\ep}\right]+\ep t\right)-\nabla u_0(z;x)\right)dtdx\right|\nonumber\\
&&\le c\int_{D^\ep}\int_Y\left|\nabla u_0\left(z;\ep\left[{x\over\ep}\right]+\ep t\right)-\nabla u_0(z;x)\right|dtdx\le c\ep.
\label{eq:a8}
\end{eqnarray}
Further,
\be
\left|\int_{D\setminus D^\ep}\int_Y\ell_i(x,y)\cdot\nabla u_0(z;x)dydx\right|
\le \int_{D\setminus D^\ep}|\nabla u_0(z;x)|dx\le  
 c\ep^{1/2}.
\label{eq:a9}
\ee
We then deduce that
\be
\left|\int_D\ell_i(x,{x\over\ep})\cdot\nabla u_0(z;x)dx-\int_D\int_Y\ell_i(x,y)\cdot\nabla u_0(z;x)dx\right|\le c\ep^{1/2}.
\label{eq:a10}
\ee
From Proposition 4.2 of \cite{HSmultirandom}, as $A(z;\cdot,\cdot)$ is uniformly bounded in $C^1(\bar D,C^1_\#(\bar Y))$, we have that $w^l$ is uniformly bounded in $C^1(\bar D,C_\#(\bar Y))$ with respect to $z\in U$.
Using 
\[
u_1(z;x,y)={\partial u_0\over\partial x_l}(z;x)w^l(z;x,y)
\]
%with $w^k$ is uniformly bounded in $C^1(\bar D, C_\#(Y))$ with respect to $z$, similarly 
we dedue that
\beqas
\left|\int_D\ell_i(x,{x\over\ep})\cdot\nabla_yu_1(z;x,{x\over\ep})dx-\int_D\int_Y\ell_i(x,y)\cdot\nabla_yu_1(z;x,y)dydx\right|\le c\ep^{1/2}.
\eeqas
This can be shown in the same way as for \eqref{eq:a10}.
 We then get the conclusion. \eproof

From Lemmas \ref{lem:1} and \ref{lem:2}, we have 
\[
|\cO_i^\ep(z)-\cO_i^0(z)|\le c\ep^{1/2},\ \ \forall\,z\in U.
\]
Thus
\[
|\cG^\ep(z)-\cG^0(z)|\le c\ep^{1/2}\ \ \forall\,z\in U.
\]
Using inequality $|\exp(-x)-\exp(-y)|\le |x-y|$ for $x,y\ge 0$, we have
\begin{eqnarray}
|\exp(-\frac12\Phi^\ep(z,\delta))-\exp(-\frac12\Phi^0(z,\delta))|\le c|\Phi^\ep(z,\delta)-\Phi^0(z,\delta)|\nonumber\\
\le c|\langle\Sigma^{-1/2}(\cG^\ep(z)-\cG^0(z)),\Sigma^{-1/2}(2\delta-\cG^\ep(z)-\cG^0(z))\rangle.
\label{eq:a11}
\end{eqnarray}
From \eqref{eq:a11}, and the uniform boundedness of $\cG^\ep(z)$ and $\cG^0(z)$ with respect to $z$, we have $I_1\le c\ep$. 
Similarly, we have $I_2\le c\ep$. Thus
\[
d_{Hell}(\rho^{\delta,\ep},\rho^{\delta,0})\le c\ep^{1/2}.
\]
\section{}\label{app:b}
We show Theorem \ref{thm:rateofconvgauss} in this appendix. For each $z\in\bar U$, the $O(\ep^{1/2})$ homogenization rate of convergence  in Lemma \ref{lem:1} holds. However, the constant $c$ now depends on $z$. We specify the dependence of this constant on $z$. 
%From \eqref{eq:tsforward}, we can write $u_1$ in terms of $u_0$. For $l=1,\ldots,d$, we denote by $w^l(z;x,y)$, as a function of $y\in Y$, the solution of the cell problem
%\be
%\nabla\cdot(A(z;x,y)\nabla w^l(z;x,y))=-\nabla\cdot (A(z;x,y)e^l),
%\label{eq:cell}
%\ee
%where $w(z;x,\cdot)\in H^1_\#(Y)$, $e^l$ is the $l$th unit vector with all the components being zero except the $l$th component which is 1. The homogenized coefficient $A^0(z;x)$ is determined by
%\be
%A^0_{kl}(z;x)=\int_YA(z;x,y)(\delta_{lp}+{\partial w^l\over\partial y_p})\delta_{kp} dy.
%\label{eq:A0}
%\ee
%The function $u_0$ is the solution of the homogenized equation 
%\be
%-\nabla\cdot(A^0(z;x)\nabla u_0(z;x))=f(x).
%\label{eq:homeq}
%\ee
%The function $u_1$ is determined by
%\be
%u_1(z;x,y)={\partial u_0\over\partial x_l}(z;x)w^l(z;x,y).
%\label{eq:u1}
%\ee

We have the following result.
\begin{lemma}\label{lem:A0}
There are constants $c_1>0$, $c_2>0$, and $c_3>0$ such that the homogenized coefficient $A^0(z;x)$ satisfies
\[
\exp(-c_1\sum_{j=1}^\infty|z_j|b_j))|\xi|^2
\le A_{kl}^0(z;x)\xi_k\xi_l\le c_2(1+\exp(c_3\sum_{j=1}^\infty|z_j|b_j))|\xi|^2
\]
for all $\xi\in\IR^d$. 
%
%\[
%|A^0_{pq}(z;x)|\le c_1(1+\exp(c_2\sum_{i=1}^\infty |z_j|b_j)).
%\]
\end{lemma}
\bproof
From \eqref{eq:A0}, we have that
\[
A^0_{kl}(z;x)\xi_k\xi_l=\int_YA(z;x,y)\left(\xi_p+{\partial (w^k\xi_k)\over\partial y_p}\right)\left(\xi_p+{\partial (w^l\xi_l)\over\partial y_p}\right)dy.
\]
We therefore have
\beqas
A^0_{kl}(z;x)\xi_k\xi_l\ge c_*(z)\int_Y\left(\xi_p+{\partial (w^k\xi_k)\over\partial y_p}\right)\left(\xi_p+{\partial (w^l\xi_l)\over\partial y_p}\right)dy\ge c_*(z)|\xi|^2.
\eeqas
On the other hand,
\beqas
A^0_{kl}(z;x)\xi_k\xi_l\le c^*(z)\left(\sum_{p=1}^d\xi_p^2+\sum_{p=1}^d\int_Y{\partial(w^l\xi_l)\over\partial y_p}{\partial (w^k\xi_k)\over\partial y_p}dy\right)\\
\le c^*(z)\left(\sum_{p=1}^d\xi_p^2+\sum_{p=1}^d\left(\sum_{k=1}^d\xi_k^2\right)\left(\sum_{k=1}^d\int_Y\left({\partial w^k\over\partial y_q}\right)^2dy\right)\right).
\eeqas
From \eqref{eq:cell}, we deduce that
\[
\|w^l(z;x,\cdot)\|_{H^1_\#(Y)/\IR}\le {c^*(z)\over c_*(z)}\|A(z;x,\cdot)e^l)\|_{L^2(Y)}\le c(1+\exp(c\sum_{j=1}^\infty|z_j|b_j)).
\]
From these we deduce
\[
A^0_{kl}(z;x)\xi_k\xi_l\le c\left(1+\exp(c\sum_{j=1}^\infty|z_j|b_j)\right)|\xi|^2.
\]
%
%From \eqref{eq:A0}, we have
%\[
%|A^0_{kl}(z;x)|\le c^*(z)(1+\|w^l(z;x,\cdot)\|_{H^1_\#(Y)/\IR})\le c(1+\exp(c\sum_{j=1}^\infty|z_j|b_j)).
%\]
\eproof

\begin{lemma}\label{lem:wlnorm}
There are constants $c_4>0$ and $c_5>0$ such that for $z\in\bar U$ 
\[
\|w^l(z)\|_{C^1(\bar D,H^2(Y))\bigcap L^\infty(\bar D,C^1(\bar Y))}\le c_4(1+\exp(c_5\sum_{i=1}^\infty |z_j|\bar b_j)).
\]
\end{lemma}
\bproof
From equation \eqref{eq:cell}, we have that
\[
\int_YA(z;x,y)\nabla_y(w^l(z;x,y)+e^l\cdot y))\cdot\nabla_y\phi dy=0,\ \ \forall\,\phi\in H^1_\#(Y).
\]
As any functions in $\cD(\IR)$ can be decomposed by partition of unity to functions of small supports which can be extended to a periodic function, this equation holds also for all $\phi\in\cD(\IR)$. Thus from theorem 4.16 of McLean \cite{McLean}, we have that 
\[
\|w^l(z;x,\cdot)+e^l\cdot y\|_{H^3(Y)}\le c\|w^l(z;x,\cdot)\|_{H^1(Y)}
\]
where $c$ only depends on the $C^{1,1}$ norm of $A$ (with respect to $y$) in a polynomial manner. It is clear that
\beqas
\|w^l(z;x,\cdot)\|_{H^1(Y)/\IR}\le {c^*(z)\over c_*(z)}\|A(z;x,\cdot)e^l\|_{L^2(Y)}\le
c(1+\exp(c\sum_{j=1}^\infty |z_j|b_j)).
\eeqas
We thus deduce that
\[
\|w^l(z;x,\cdot)\|_{H^3(Y)/\IR}\le c(1+\exp(c\sum_{j=1}^\infty |z_j|\bar b_j)).
\]
For $d=2,3$, $H^3(Y)\subset C^1(\bar Y)$ so $\|w^l(z;x,\cdot)\|_{C^1(\bar Y)}\le c(1+\exp(c\sum_{j=1}^\infty|z_j|\bar b_j))$.
For $l,k=1,\ldots,d$
\[
-\nabla_y\cdot\left(A\nabla_y{\partial w^l\over\partial x_k}\right)=\nabla_y\cdot\left({\partial A\over\partial x_k}(e^l+\nabla_y w^l)\right).
\]
Therefore,
\[
\left\|{\partial w^l\over\partial x_k}\right\|_{H^1(Y)/\IR}\le {c^*(z)\over c_*(z)}\left\|{\partial A\over\partial x_k}(e^l+\nabla w^l)\right\|_{L^2(Y)}\le c(1+\exp(c\sum_{j=1}^\infty|z_j|\bar b_j)).
\]
Thus from Theorem 4.16 of \cite{McLean}, we have 
\[
\left\|{\partial w^l\over\partial x_k}(z;x,\cdot)\right\|_{H^2(Y)/\IR}\le c\left(\|{\partial w^l\over\partial x_k}(z;x,\cdot)\|_{H^1(Y)/\IR}+\left\|{\partial A\over\partial x_k}(e^l+\nabla w^l(z;x,\cdot))\right\|_{H^1(Y)}\right)
\]
where the constant $c$ depends only on $c^*(z)$, $c_*(z)$ and the $C^{0,1}$ norm of $A(z;x,\cdot)$ polynomially. We thus deduce
\[
\left\|{\partial w^l\over\partial x_k}(z;x,\cdot)\right\|_{H^2(Y)/\IR}\le c\left(1+\exp(c\sum_{j=1}^\infty |z_j|\bar b_j)\right).
\]
\eproof

\begin{lemma}\label{lem:H2normu0}
Assume that $D$ is a convex domain, and $f\in L^2(D)$. Under Assumption \ref{assum:furtherregularitypsigauss}, there are constants $c_6>0$ and $c_7>0$ such that 
for all $z\in\bar U$
\[
\|u_0(z;\cdot)\|_{H^2(D)}\le c_6\left(1+\exp(c_7\sum_{i=1}^\infty |z_j|\bar b_j)\right).
\]
\end{lemma}
\bproof
The solution $u_0(z;\cdot)$ of \eqref{eq:homeq} belongs to $H^2(D)$ when the domain $D$ is convex. Examining the proof of Theorem 3.1.3.1 of Grisvard \cite{Grisvard}, we find that  
\[
\|u_0(z;\cdot)\|_{H^2(D)}\le c\|f\|_{L^2(D)}
\]
where the constant $c$ depends polynomially on the Lipschitz norm of $A^0_{kl}(z;\cdot)$, the upper bound of the entries of $R=A^0(z;x)^{-1/2}$, and the constant of the Friedrichs inequality in the domain $RD$. From \eqref{eq:A0} and Lemma \ref{lem:wlnorm}, the Lipschitz norm of $A^0_{kl}$ is bounded by
$
c(1+\exp(c\sum_{j=1}^\infty|z_j|\bar b_j)).
$
From Lemma \ref{lem:A0}, the eigenvalues of $A^0$ are bounded below by $c\exp(-c\sum_{j=1}^\infty|z_j|b_j))$ and bounded above by $c(1+\exp(c\sum_{j=1}^\infty|z_j|b_j))$. Thus the eigenvalues of $A^0(z;x)^{-1/2}$ is bounded %below by $\exp(-c\sum_{j=1}^\infty|z_j|b_j)$ and bounded 
above by $c\exp(c\sum_{j=1}^\infty|z_j|b_j)$. Let $\xi$ be the vector in $\IR^d$ whose components are all zero except the $k$th and the $l$th components which are 1. As $A^0(z;x)^{-1/2}$ is symmetric,
\[
(A^0)^{-1/2}_{kl}(z;x)={1\over 2}({A^0})^{-1/2}(z;x)\xi\cdot\xi
\]
so the entries of $(A^0)^{-1/2}(z;x)$ are bounded above by $c(1+\exp(c\sum_{j=1}^\infty|z_j|b_j))$. We note that the constant in the Friedrichs inequality in $RD$ is bounded polynomially by the diameter of $RD$ (see, e.g., Wloka \cite{Wloka} page 116) so is also bounded by $c(1+\exp(c\sum_{j=1}^\infty|z_j|b_j))$. We get the conclusion.\eproof

\begin{proposition}\label{prop:homerrorgauss}
Assume that the domain $D$ is convex, and $f\in L^2(D)$. Under Assumption \ref{assum:furtherregularitypsigauss}, there are constants $c_8>0$ and $c_9>0$ such that for $z\in\bar U$
\[
\|\nabla\ue(z)-[\nabla u_0(z)+\nabla_y u_1(z;\cdot,{\cdot\over\ep})\|_{L^2(D)}\le C(z)\ep^{1/2},
\]
where
\[
C(z)=c_8(1+\exp(c_9\sum_{i=1}^\infty |z_j|\bar b_j)).
\]
\end{proposition}
\bproof
We check the dependence  on the parameters of the constants in the proof of the homogenization convergence rate for parametric two scale elliptic problem. We follow closely the proof of Proposition 5.1 in \cite{HSmultirandom}. The proof is an extension of the homogenization proof in \cite{JKO} for the case of smooth $u_0$ to the case where $u_0$ is only in $H^2(D)$. 

 Let
\[
u_1^\ep(z;x)=u_0(z;x)+\ep w^l(z;x,{x\over\ep}){\partial u_0(z;x)\over\partial x_l}.
\]
We first show that
\[
\|{\rm div}A^\ep\nabla u_1^\ep-{\rm div}A^0\nabla u_0\|_{H^{-1}(D)}\le c(z)\ep,
\]
where $c(z)$ is of the form 
\be
c(1+\exp(c\sum_{j=1}^\infty|z_j|\bar b_j)).
\label{eq:formofc}
\ee
We note that
\beqas
(A^\ep(z;x)\nabla u_1^\ep(z;x))_i=\Bigl(A^\ep_{ij}(z;x)+A^\ep_{ik}(z;x){\partial w^j\over\partial y_k}(z;x,{x\over\ep})\Bigr){\partial u_0\over\partial x_j}(z;x)+\ep A^\ep_{ij}(z;x)w^k(z;x,{x\over\ep}){\partial^2u_0\over\partial x_j\partial x_k}(z;x)\\
=A^0_{ij}(z;x){\partial u_0\over\partial x_j}(z;x)+\Bigl(A^\ep_{ij}(z;x)+A^\ep_{ik}(z;x){\partial w^j\over\partial y_k}(z;x,{x\over\ep})-A^0_{ij}(z;x)\Bigr){\partial u_0\over\partial x_j}(z;x)\\
+\ep A^\ep_{ij}(z;x)w^k(z;x,{x\over\ep}){\partial^2u_0\over\partial x_j\partial x_k}(z;x)\\
=A^0_{ij}(z;x){\partial u_0\over\partial x_j}(z;x)+g_i^j(z;x,{x\over\ep}){\partial u_0\over\partial x_j}(z;x)+\ep A^\ep_{ij}(z;x)w^k(z;x,{x\over\ep}){\partial^2u_0\over\partial x_j\partial x_k}(z;x),
\eeqas
where the functions $g_i^j(z;x,y)$ which are periodic in $y$ are defined as
\[
g_i^j(z;x,y)=A_{ij}(z;x,y)+A_{ik}(z;x,y){\partial w^j\over\partial y_k}(z;x,y)-A^0_{ij}(z;x).
\]
We have
\[
\int_Yg_i^j(z;x,y)dy=0,
\mbox{ and }
{\partial\over\partial y_i}g_i^j(z;x,y)=0.
\]
Followng Jikov et al. \cite{JKO}, we write $g_i^k(z;x,y)={\partial\over\partial y_j}\alpha_{ij}^k(z;x,y)$ where $\alpha_{ij}^k(z;x,y)$ are periodic in $y$ and $\alpha_{ij}^k=-\alpha_{ji}^k$. As $w^l(z;\cdot,\cdot)\in C^1(\bar D,H^2_\#(Y))$ and $A(z;\cdot,\cdot)\in C^1(\bar D,C^1_\#(\bar Y))$, $g_i^j(z;\cdot,\cdot)\in C^1(\bar D,H^1_\#(Y))$. From Lemma \ref{lem:wlnorm}, the $C^1(\bar D, H^1(\bar Y))$ norm of $g_i^j(z;\cdot,\cdot)$ is bounded above a constant of the form
\eqref{eq:formofc}
% Jikov et al (\cite{JKO} page 7) construct the functions $\alpha_{ij}^k$ as follows. 
Consider the  Fourier series of $\bg^k$
\[
\bg^k(z;x,y)=\sum_{l\in\spa,l\ne 0}\bg^k_l(z;x)\exp(\sqrt{-1}l\cdot y).
\]
As $\bg^k(z;\cdot,\cdot)\in C^1(\bar D,H^1_\#(Y))^d$, $\bg^k_l(z;\cdot)\in C^1(\bar D)^d$ and, for all $r=1,\ldots,d$, there is a constant $c_{r}(z)$ such that 
\be
\sum_{l\in\spa,l\ne 0}|\bg^k_l(z;x)|^2l_r^2<c_{r}(z)
\label{eq:b1}
\ee
where $c_{r}(z)$
%<c\|\bg^k\|_{C^1(\bar D,H^1_\#(Y))^d}$ 
is of the form \eqref{eq:formofc} due to the bound of  $\|\bg^k(z;\cdot,\cdot)\|_{C^1(\bar D,H^1_\#(Y))^d}$. From \cite{JKO}, the functions $\alpha_{ij}^k$ are defined as
\[
\alpha^k_{ij}(z;x,y)=-\sqrt{-1}\sum_{l\in\spa,l\ne 0}{(\bg^k_l(z;x))_jl_i-(\bg^k_l(z;x))_il_j\over|l|^2}\exp(\sqrt{-1}l\cdot y).
\]
%It is then obvious that f
For $r,s=1,\ldots,d$
\[
\sum_{l\in\spa,l\ne 0}{|(\bg^k_l(z;x))_jl_i-(\bg^k_l(z;x))_il_j|^2\over|l|^4}l_r^2l_s^2<c_{rs}(z).
\]
From \eqref{eq:b1}, the constant $c_{rs}(z)$ can be taken to be of the form \eqref{eq:formofc}. 
Therefore for dimension $d\le 3$, $\alpha^k_{ij}(z;\cdot,\cdot)\in C^1(\bar D,H^2_\#(Y)\subset C^1(\bar D,C(\bar Y))$; and the $C^1(\bar D, C(\bar Y))$ norm of $\alpha^k_{ij}(z;\cdot,\cdot)$ is bounded above by a constant of the form \eqref{eq:formofc}. 

We note that
\[
(A^\ep\nabla u_1^\ep(z;x)-A_0\nabla u_0(z;x))_i=\ep{\partial\over\partial x_j}\Bigl(\alpha^k_{ij}(z;x,{x\over\ep}){\partial u_0(z;x)\over\partial x_k}\Bigr)+(r_\ep)_i(z;x),
\]
where
\[
(r_\ep)_i(z;x)=-\ep{\partial\alpha^k_{ij}(z;x,y)\over\partial x_j}\Big|_{y=x/\ep}{\partial u_0(z;x)\over\partial x_k}-\ep\alpha_{ij}^k(z;x,{x\over\ep}){\partial^2 u_0(z;x)\over\partial x_k\partial x_j}+\ep w^k(z;x,{x\over\ep})A^\ep_{ij}(z;x,{x\over\ep}){\partial^2u_0(z;x)\over\partial x_k\partial x_j}.
\]
As $\alpha^k_{ij}(z;\cdot,\cdot)\in C^1(\bar D,C(\bar Y))$, $\|(r_\ep)_i(z;\cdot)\|_{L^2(D)}\le c(z)\ep$ where the constant $c(z)$ is of the form \eqref{eq:formofc}. As $\alpha^k_{ij}=\alpha^k_{ji}$,
\[
\|{\rm div}A^\ep\nabla u_1^\ep(z;x)-{\rm div}A^0\nabla u_0(z;x)\|_{H^{-1}(D)}\le c(z)\ep,
\]
where the constant $c(z)$ is of the form \eqref{eq:formofc}. 

As ${\rm div}A^0\nabla u_0={\rm div}A^\ep\nabla u^\ep$, we find that
\[
\|{\rm div}A^\ep\nabla u_1^\ep(z;x)-{\rm div}A^\ep\nabla u^\ep(z;x)\|_{H^{-1}(D)}\le c(z)\ep,
\] 
where the constant $c(z)$ is of the form \eqref{eq:formofc}.

Let $\tau^\ep\in{\cal D}(D)$ be such that $\tau^\ep=1$ outside an $\ep$ neighbourbood of $\partial D$ and $\ep|\nabla\tau^\ep(x)|\le c$ for all $\ep$. We consider the function
\[
w^\ep_1(z;x)=u_0(z;x)+\ep\tau^\ep(x)w^k(z;x,{x\over\ep}){\partial u_0(z;x)\over\partial x_k}=u_1^\ep(z;x)-\ep(1-\tau^\ep(x))w^k(z;x,{x\over\ep}){\partial u_0(z;x)\over\partial x_k}.
\]
We then get
\beqas
{\partial\over\partial x_j}(u_1^\ep-w_1^\ep)(z;x)=-\ep{\partial\tau^\ep\over\partial x_j}(x)w^k(z;x,{x\over\ep}){\partial u_0(z;x)\over\partial x_k}+(1-\tau^\ep(x)){\partial w^k\over\partial y_j}(z;x,{x\over\ep}){\partial u_0\over\partial x_k}(z;x)+\\
\ep(1-\tau^\ep(x))w^k(z;x,{x\over\ep}){\partial^2u_0(z;x)\over\partial x_k\partial x_j}.
\eeqas
Let $D^\ep$ be the $\ep$ neighbourhood of $\partial D$. As $\partial D$ is Lipschitz; for all smooth functions $\phi(x)$
\[
\|\phi\|^2_{L^2(D^\ep)}\le c\ep^2\|\phi\|^2_{H^1(D)}+c\ep\|\phi\|^2_{L^2(\partial D)},
\]
where $c$ only depends on the domain $D$; so for all $\phi\in H^1(D)$
\[
\|\phi\|_{L^2(D^\ep)}\le c\ep^{1/2}\|\phi\|_{H^1(D)}.
\]
Therefore from Lemmas \ref{lem:H2normu0} and \ref{lem:wlnorm} we get 
%since $u_0(z;x)\in L^\infty(U,H^2(D))$ and $w^l(z;x,y)\in L^\infty(U\times\bar D,C^1(\bar Y))$
\[
\|u_1^\ep(z;\cdot)-w_1^\ep(z;\cdot)\|_{H^1(D)}\le c(z)\ep^{1/2},
\]
where the constant $c(z)$ is of the form \eqref{eq:formofc} . Thus,
\[
\|{\rm div}(A^\ep(\nabla u_1^\ep(z;\cdot)-\nabla w_1^\ep(z;\cdot)))\|_{H^{-1}(D)}\le c(z)\ep^{1/2},
\]
so
\[
\|{\rm div}(A^\ep(\nabla u^\ep(z;\cdot)-\nabla w_1^\ep(z;\cdot)))\|_{H^{-1}(D)}\le c(z)\ep^{1/2}
\]
where the constant $c(z)$ is of the form \eqref{eq:formofc}. 
From \eqref{eq:coercivebounded}, we get
\[
\|u^\ep(z;\cdot)-w_1^\ep(z;\cdot)\|_{H^1_0(D)}\le c(z)\ep^{1/2},
\]
where $c(z)$ is of the form \eqref{eq:formofc}. Hence
\[
\|u^\ep(z;\cdot)-u_1^\ep(z;\cdot)\|_{H^1(D)}\le c(z)\ep^{1/2}
\]
where $c(z)$ is of the form \eqref{eq:formofc}. 
\eproof

{\it Proof of Theorem \ref{thm:rateofconvgauss}}

The proof of Theorem \ref{thm:rateofconvgauss} is largely similar to that in Appendix \ref{app:a} except that the constants depend on $z$. 
%We will show that these constants are bounded by 
%\be
%c(1+\exp(c\sum_{i=1}^\infty |z_j|\bar b_j))
%\label{eq:formofc}
%\ee
%for $z\in\bar U$. 

Consider $z\in\bar U$. From Proposition \ref{prop:homerrorgauss}, the constant in \eqref{eq:a1} is of the form \eqref{eq:formofc}. From Lemma \ref{lem:H2normu0}, the constant $c$ in \eqref{eq:a3} is also of the form \eqref{eq:formofc}. The proof of Lemma 5.5 of \cite{HSmultirandom} shows that the constant $c$ in \eqref{eq:a4} is of the form $c\|u_0\|_{H^2(D)}$ so is of the form \eqref{eq:formofc}. The constant $c$ in \eqref{eq:a6} depends on $\int_D\int_Y\nabla u_0(z;\ep[{x\over\ep}]+\ep t)dtdx$. From \eqref{eq:a4}, this is bounded by $\|\nabla u_0\|_{L^2(D)}$ so the constant $c$ in \eqref{eq:a6} is of the form \eqref{eq:formofc}. Similarly, the constant $c$ in equations \eqref{eq:a7}, \eqref{eq:a8},\eqref{eq:a9}, and \eqref{eq:a10} are all of the form \eqref{eq:formofc}. Arguing in the same way, from Lemma \ref{lem:wlnorm}, we have that 
\beqas
\left|\int_D\ell_i(x,{x\over\ep})\cdot\nabla_yu_1(z;x,{x\over\ep})dx-\int_D\int_Y\ell_i(x,y)\cdot\nabla_yu_1(z;x,y)dydx\right|\le c(z)\ep^{1/2}
\eeqas
where $c(z)$ is of the form \eqref{eq:formofc}. Thus for $z\in \bar U$, Lemma \ref{lem:2} holds for the case of Gaussian prior with the constant $c$ depending on $z$ and is of the form \eqref{eq:formofc}. 

We therefore deduce that 
\be
|\cG^\ep(z)-\cG^0(z)|\le c(z)\ep^{1/2}
\label{eq:bb}
\ee
where $c(z)$ is of the form \eqref{eq:formofc}. From \eqref{eq:bounduegauss}, we have that 
\be
|\cG^\ep(z)|\le c(z)
\label{eq:bc}
\ee
where $c(z)$ is of the form \eqref{eq:formofc}. From Lemmas %\ref{lem:H2normu0} and  
\ref{lem:wlnorm} and \ref{lem:H2normu0}, and equation \eqref{eq:u1}, we deduce that
\be
|\cG^0(z)|\le c(z)
\label{eq:bd}
\ee
where $c(z)$ is of the form \eqref{eq:formofc}. %Using inequality $|\exp(-x)-\exp(-y)|\le |x-y|$ for $x,y\ge 0$, we have
%\begin{eqnarray}
%|\exp(-\frac12\Phi^\ep(z,\delta))-\exp(-\frac12\Phi^0(z,\delta))|\le c|\Phi^\ep(z,\delta)-\Phi^0(z,\delta)|\nonumber\\
%\le c|\langle\Sigma^{-1/2}(\cG^\ep(z)-\cG^0(z)),\Sigma^{-1/2}(2\delta-\cG^\ep(z)-\cG^0(z))\rangle.
%\label{eq:11}
%\end{eqnarray}
From \eqref{eq:bb}, \eqref{eq:bc} and \eqref{eq:bd}, the right hand side of \eqref{eq:a11} is bounded by $c(z)\ep^{1/2}\ \forall\,z\in\bar U$ where $c(z)$ is of the form \eqref{eq:formofc}. Thus from \eqref{eq:I1}, as the set $\bar U$ has $\rho$ measure 1 the integral in \eqref{eq:I1} over $U\setminus\bar U$ is 0, we deduce that $I_1\le c\ep$ (see Lemma \ref{lem:B100} below). Similarly, $I_2\le c\ep$. The conclusion then follows. 
\begin{lemma}\label{lem:B100}
For any constants $\alpha$ and $\beta$, 
\[
\int_{\bar U} \alpha(1+\exp(\beta\sum_{j=1}^\infty|z_j|\bar b_j))d\rho(z)
\]
is finite. 
\end{lemma}
\bproof We use the inequality for $s>0$
\[
\int_{-\infty}^\infty\exp(-t^2/2+|t|s){dt\over\sqrt{2\pi}}\le\exp(s^2/2)\exp(s\sqrt{2/\pi})
\]
which is proved in \cite{HSmcmclogn}. We have
\beqas
\int_{\bar U}\alpha (1+\exp(\beta\sum_{j=1}^\infty|z_j|\bar b_j))d\rho(z)\le\alpha+\alpha\exp\left(\frac12\beta^2\sum_{j=1}^\infty\bar b_j^2\right)\exp\left(\sqrt{2\over\pi}\beta\sum_{j=1}^\infty\bar b_j\right)
\eeqas
which is finite.
\eproof

\end{appendices}

\bibliographystyle{plain}
\bibliography{inv_hom}

\end{document}